# ERROR ESTIMATES FOR BINOMIAL APPROXIMATIONS OF GAME OPTIONS

By Yuri Kifer

*Hebrew University*

We justify and give error estimates for binomial approximations of game (Israeli) options in the Black–Scholes market with Lipschitz continuous path dependent payoffs which are new also for usual American style options. We show also that rational (optimal) exercise times and hedging self-financing portfolios of binomial approximations yield for game options in the Black–Scholes market "nearly" rational exercise times and "nearly" hedging self-financing portfolios with small average shortfalls and initial capitals close to fair prices of the options. The estimates rely on strong invariance principle type approximations via the Skorokhod embedding.

**1. Introduction.** Cox, Ross and Rubinstein's (CRR) binomial model of a financial market was introduced in [6] not only as a simplified discrete time and space counterpart of the Black–Scholes (BS) model based on the geometric Brownian motion, but also as a convenient approximation of the latter which, indeed, became a popular tool to evaluate various derivative securities. Clearly, for an approximation to have any practical value, it is necessary to estimate the corresponding error. Many papers dealt with both justification of the CRR approximation of European and American options in the BS market and with estimates of the corresponding error terms (see, e.g., [1, 10, 21, 22, 24, 27, 30]). Still, none of these papers derived error estimates for options with path dependent payoffs. In this paper we consider game (Israeli) options introduced in [15] which generalize American style options when not only their holders but also their writers have the right to exercise and we obtain error estimates of binomial approximations of fair prices, rational exercise times and hedging self-financing portfolios for such options considered in a BS market and having path dependent payoffs.

A game option (or contingent claim) studied in [15] is a contract between a writer and a holder at time $t = 0$ such that both have the right to exercise at any stopping time before the expiry date $T$. If the holder exercises at time $t$, he may claim the amount $Y_t \geq 0$ from the writer and if the writer exercises at time $t$, he must pay to the holder the amount $X_t \geq Y_t$ so that $\delta_t = X_t - Y_t$ is viewed as a penalty imposed on the writer for cancellation of the contract. If both exercise at the same time $t$, then the holder may claim $Y_t$ and if neither have exercised until the expiry time $T$, then the holder may claim the amount $Y_T$. In short, if the writer will exercise at a stopping time $\sigma \leq T$ and the holder at a stopping time $\tau \leq T$, then the former pays to the latter the amount $R(\sigma, \tau)$, where

$$R(s,t) = X_s \mathbb{I}_{s<t} + Y_t \mathbb{I}_{t \leq s} \tag{1.1}$$

and we set $\mathbb{I}_A = 1$ if (an event or an assertion) $A$ holds true and $\mathbb{I}_A = 0$ if not. As usual, we start with a complete probability space $(\Omega, \mathcal{F}, P)$ and a filtration of $\sigma$-algebras $\{\mathcal{F}_t\}_{t\geq 0}$ generated either by a Brownian motion in the BS model or by i.i.d. binomial random variables in the CRR model. The payoff processes $X_t$ and $Y_t$ should be adapted to the corresponding filtration and in the continuous time case they are supposed to be right continuous with left limits, though the latter could be relaxed sometimes.

Two popular models of complete markets were considered in [15] for pricing of game options. First, the discrete time CRR binomial model was treated there where the stock price $S_k$ at time $k$ is equal to

$$S_k = S_0 \prod_{j=1}^{k} (1+\rho_j), \qquad S_0 > 0, \tag{1.2}$$

where $\rho_j$, $j = 1, 2, \ldots$, are independent identically distributed (i.i.d.) random variables such that $\rho_j = b > 0$ with probability $p > 0$ and $\rho_j = a < 0, a > -1$ with probability $q = 1 - p > 0$. Second, [15] deals with the continuous time BS market model where the stock price $S_t$ at time $t$ is given by the geometric Brownian motion

$$S_t = S_0 \exp((\alpha - \kappa^2/2)t + \kappa B_t), \qquad S_0 > 0, \tag{1.3}$$

where $\{B_t\}_{t\geq 0}$ is the standard one-dimensional continuous in time Brownian motion (Wiener process) starting at zero and $\kappa > 0$, $\alpha \in (-\infty, \infty)$ are some parameters. In addition to the stock which is a risky security, the market includes in both cases also a savings account with a deterministic growth given by the formulas

$$b_n = (1+r)^n b_0 \quad \text{and} \quad b_t = b_0 e^{rt}, \qquad b_0, r > 0, \tag{1.4}$$

in the CRR model (where we assume, in addition, that $r < b$) and in the BS model, respectively.

Recall (see [28]) that a probability measure describing the evolution of a stock price in a stochastic financial market is called martingale (risk-neutral) if the discounted stock prices $[(1+r)^{-k}S_k$ in the CRR model and $e^{-rt}S_t$ in the BS model] become martingales. Relying on hedging arguments, it was shown in [15] that the fair price $V$ of the game option is given by the formulas

$$V = \min_{\sigma \in \mathcal{T}_{0T}} \max_{\tau \in \mathcal{T}_{0T}} E((1+r)^{-\sigma\wedge\tau} R(\sigma,\tau)) \tag{1.5}$$

in the CRR market [with usual notation $a \wedge b = \min(a,b)$, $a \vee b = \max(a,b)$] and

$$V = \inf_{\sigma \in \mathcal{T}_{0T}} \sup_{\tau \in \mathcal{T}_{0T}} E(e^{-r\sigma\wedge\tau} R(\sigma,\tau)) \tag{1.6}$$

in the BS market, where the expectations are taken with respect to the corresponding martingale probabilities, which are uniquely defined since these markets are known to be complete (see [28]), $T$ is the expiry time and $\mathcal{T}_{st}$ is the space of corresponding stopping times with values between $s$ and $t$ taking into account that in the CRR model $\sigma$ and $\tau$ are allowed to take only integer values. Observe that formulas (1.5) and (1.6) represent also the values of corresponding Dynkin's (optimal stopping) games with payoffs $(1+r)^{-\sigma\wedge\tau}R(\sigma,\tau)$ and $e^{-r\sigma\wedge\tau}R(\sigma,\tau)$, respectively, when the first and the second players stop the game at stopping times $\sigma$ and $\tau$, respectively. Observe that since their introduction in [15], various aspects of game (Israeli) options were studied in [2, 5, 7, 9, 12, 17, 18, 19, 20] and recently this technique was applied in [8] to convertible (callable) bonds.

The continuous time BS model is generally considered as a better description of the evolution of real stocks, in particular, since the CRR model allows only two possible values $(1+b)S_k$ and $(1+a)S_k$ for the stock price $S_{k+1}$ at time $k+1$ given its price $S_k$ at time $k$. The main advantage of the CRR model is its simplicity and the possibility of easier computations of the value $V$ in (1.5), in particular, by means of the dynamical programming recursive relations (see [15]),

$$\begin{aligned} &V = V_{0,N}, \qquad V_{N,N} = (1+r)^{-N} Y_N \quad \text{and} \\ &V_{k,N} = \min((1+r)^{-k} X_k, \max((1+r)^{-k} Y_k, E(V_{k+1,N} | \mathcal{F}_k))), \end{aligned} \tag{1.7}$$

where a positive integer $N$ is an expiry time and $\{\mathcal{F}_k\}_{k\ge 0}$ is the corresponding filtration of $\sigma$-algebras.

Following [30], we will approximate the BS model by a sequence of CRR models with the interest rates $r = r^{(n)}$ from (1.4) and with random variables $\rho_k = \rho_k^{(n)}$ from (1.2) given by

$$\begin{aligned} &r = r^{(n)} = \exp(rT/n) - 1 \quad \text{and} \\ &\rho_k = \rho_k^{(n)} = \exp\left(\frac{rT}{n} + \kappa\left(\frac{T}{n}\right)^{1/2} \xi_k\right) - 1, \end{aligned} \tag{1.8}$$

where $\xi_j = \xi_j^{(n)}, j = 1, 2, \ldots,$ are i.i.d. random variables taking on the values 1 and $-1$ with probabilities $p^{(n)} = (\exp(\kappa\sqrt{\frac{T}{n}}) + 1)^{-1}$ and $1 - p^{(n)} = (\exp(-\kappa\sqrt{\frac{T}{n}}) + 1)^{-1}$, respectively. This choice of random variables $\xi_i$, $i \in \mathbb{N}$, determines already the probability measures $P_n^\xi = \{p^{(n)}, 1 - p^{(n)}\}^\infty$ for the above sequence of CRR models and since $E_n^\xi \rho_k^{(n)} = r^{(n)}$, where $E_n^\xi$ is the expectation with respect to $P_n^\xi$, we conclude that $P_n^\xi$ is the martingale measure for the corresponding CRR market and the fair price $V = V^{(n)}$ of a game option in this market is given by the formula (1.5) with $E = E_n^\xi$. Some authors consider a bit simpler and more straightforward approximation (see, e.g., [21] and [22]) where

$$\rho_k = \hat{\rho}_k = \exp\left(\frac{\beta T}{n} + \kappa\left(\frac{T}{n}\right)^{1/2} \hat{\xi}_k\right) - 1 \tag{1.9}$$

with $\beta = r - \kappa^2/2$ and $\hat{\xi}_k$ takes on the values 1 and $-1$ with the same probability 1/2. This approximation leads to similar errors estimates (with, essentially, the same proof) but, in general, we do not arrive at the martingale probability measures in this case. Thus, we have to speak then about the prices $V^{(n)}$ of discrete time Dynkin's games (rather than about the fair prices of the corresponding game options) given by (1.5) for the CRR market with $\hat{\rho}_k$ and $r^{(n)}$ described above, but with the expectation $E = E^{\hat{\xi}}$ taken with respect to the probability $P^{\hat{\xi}}$ generated by $\hat{\xi}_j$, $j = 1, 2, \ldots,$ rather than with respect to the corresponding martingale probabilities. For purposes of approximation, this difference is not so important, but the first approximation becomes more convenient for the construction of self-financing "nearly" hedging portfolios with small average shortfalls. Another useful advantage of the first approximation is that it leads to discounted stock prices evolving on the multiplicative lattice $\{S_0 \exp(m\kappa(T/n)^{1/2}), m \in \mathbb{Z}\}$ which substantially simplifies computations.

Let $V$ be the fair price of the game option in the BS market. The main goal of this paper is to show that for a certain natural class of payoffs $X_t$ and $Y_t$ which may depend on the whole path (history) of the stock price evolution (as in integral or Russian type options) the error $|V - V^{(n)}|$ does not exceed $Cn^{-1/4}(\ln n)^{3/4}$, where $C > 0$ does not depend on $n$ and it can be estimated explicitly. Moreover, we will show that the rational exercise times of our CRR binomial approximations yield near rational [$(Cn^{-1/4}(\ln n)^{3/4})$-optimal stopping times for the corresponding Dynkin games] exercise times for game options in the BS market. Since the values $V^{(n)}$ and the optimal stopping times of the corresponding discrete time Dynkin's games can be obtained directly via the dynamical programming recursive procedure (1.7), our results provide a justification of a rather effective method of computation of fair prices and exercise times of game options with path dependent payoffs. The standard construction of a self-financing hedging portfolio involves usually the Doob–Meyer decomposition of supermartingales which is explicit only in the discrete, but not in the continuous time case. We will show how to construct a self-financing portfolio in the BS market with a small average (maximal) shortfall and an initial capital close to the fair price of a game option using hedging self-financing portfolios for the approximating binomial CRR markets. The latter problem does not seem to have been addressed until now in the literature on this subject. This hints, in particular, that since hedging self-financing portfolio strategies can be computed only approximately, their possible shortfalls come naturally into the picture and they should be taken into account in option pricing even if a perfect hedging exists theoretically. Note that the results of the present paper require not only an approximation of stock prices and the corresponding payoffs, but also we have to take care about the different nature of stopping times in (1.5) and (1.6). It would be interesting to obtain similar results for discrete time and space (say, multinomial) approximations of sufficiently general Lévy markets, that is, markets where the stock price evolve according to a geometric Lévy process with jumps, but this requires additional ideas and machinery. Some discrete time approximation results without error estimates for American options in the Merton stock market model were obtained in [25]. Game options with jump–diffusion models of stock evolutions were considered recently in [8].

Our main tool is the Skorokhod type embedding of sums of i.i.d. random variables into a Brownian motion (with a constant drift, in our case). This tool was already employed for similar purposes in [24] and [30]. The first paper treats an optimal stopping problem which can be applied to an American style option with a payoff function depending only on the current stock price and, more importantly, this function must be bounded and have two bounded derivatives which excludes usual put and call options. The second paper deals only with European options and, again, only payoffs (though with some discontinuities) determined by the current stock price are allowed. A number of other papers also deal with error estimates for the CRR approximation of European and American option prices in the BS market (see, e.g., [21, 22, 27] and references there), but none of them treat path dependent payoffs (moreover, boundedness of payoffs conditions there usually exclude even American style call options) and none of them consider the game options case as well.

The main results of this paper are formulated in the next section where we discuss also the Skorokhod type embedding which we employ in the proof. In Section 3 we show how this embedding enables us to consider both CRR and BS stock evolutions in an appropriate way on the same probability space and we exhibit there a series of steps which lead to the proof of the main results. The necessary technical estimates are derived in Section 4. In Section 5 we deal with rational exercise times and self-financing nearly hedging portfolios with small averaged shortfalls. In Section 6 we generalize to the game options situation the estimates from [24], which cannot be applied to the standard options as the proof relies on very restrictive bondedness and smoothness assumptions, but still, in view of their simplicity, the arguments there may have a pedagogical value and some readers may prefer to read this case first.

**2. Preliminaries and main results.** For each $t>0$, denote by $M[0,t]$ the space of Borel measurable functions on $[0,t]$ with the uniform metric $d_{0t}(\upsilon,\tilde{\upsilon})=\sup_{0\le s\le t}|\upsilon_s-\tilde{\upsilon}_s|$. For each $t>0$, let $F_t$ and $\Delta_t$ be nonnegative functions on $M[0,t]$ such that, for some constant $L\ge 1$ and for any $t\ge s\ge 0$ and $\upsilon,\tilde{\upsilon}\in M[0,t]$,

$$|F_s(\upsilon)-F_s(\tilde{\upsilon})|+|\Delta_s(\upsilon)-\Delta_s(\tilde{\upsilon})|\le L(s+1)\,d_{0s}(\upsilon,\tilde{\upsilon}) \tag{2.1}$$

and

$$\begin{aligned}&|F_t(\upsilon)-F_s(\upsilon)|+|\Delta_t(\upsilon)-\Delta_s(\upsilon)|\\ &\qquad\le L\biggl(|t-s|\biggl(1+\sup_{u\in[0,t]}|\upsilon_u|\biggr)+\sup_{u\in[s,t]}|\upsilon_u-\upsilon_s|\biggr).\end{aligned} \tag{2.2}$$

By (2.1), $F_0(\upsilon)=F_0(\upsilon_0)$ and $\Delta_0(\upsilon)=\Delta_0(\upsilon_0)$ are functions of $\upsilon_0$ only. By (2.2),

$$F_t(\upsilon)+\Delta_t(\upsilon)\le F_0(\upsilon_0)+\Delta_0(\upsilon_0)+L(t+2)\biggl(1+\sup_{0\le s\le t}|\upsilon_s|\biggr). \tag{2.3}$$

Next, we consider the BS market on a complete probability space together with its martingale measure $P^B$ which exists and is unique as a corollary of the Girsanov theorem (see [28]). Let $B_t$, $t\ge 0$, be the standard one-dimensional continuous in time Brownian motion with respect to the martingale measure $P^B$. Set

$$B_t^*=-\frac{\kappa}{2}t+B_t,\qquad t\ge 0.$$

Then the stock price $S_t^B(z)$ at time $t$ in the BS market can be written in the form

$$S_t^B(z)=z\exp(rt+\kappa B_t^*),\qquad S_0^B(z)=z>0, \tag{2.4}$$

where $r>0$ is the interest rate and $\kappa>0$ is the so-called volatility. We will consider game options in the BS market with payoff processes in the form

$$Y_t=F_t(S^B(z))\quad\text{ and }\quad X_t=G_t(S^B(z)),\qquad t\in[0,T],T>0,$$

where $G_t=F_t+\Delta_t$, $F,\Delta$ satisfy (2.1) and (2.2), $S^B(z)=S^B(z,\omega)\in M[0,T]$ is a random function taking the value $S_t^B(z)=S_t^B(z,\omega)$ at $t\in[0,T]$, and in the notation $F_t(S^B(z))$, $G_t(S^B(z))$ for $t<T$, we take the restriction of $S^B(z)$ to the interval $[0,t]$. The fair price $V=V(z)$ of this option with an initial value $z>0$ of the stock is given by (1.6).

Next, we consider a sequence of CRR markets on a complete probability space such that, for each $n=1,2,\dots$, the stock prices $S_t^{(n)}(z)$ at time $t$ are given by the formula

$$
\begin{aligned}
&S_t^{(n)}(z)=z\exp\left(\sum_{k=1}^{[nt/T]}\left(\frac{rT}{n}+\kappa\left(\frac{T}{n}\right)^{1/2}\xi_k\right)\right),\qquad t\ge T/n\quad\text{and}\\
&S_t^{(n)}(z)=S_0^{(n)}(z)=z>0,\qquad t\in[0,T/n)
\end{aligned}
\tag{2.5}
$$

where, recall, $\xi_1,\xi_2,\dots$ are i.i.d. random variables taking the values 1 and $-1$ with probabilities $p^{(n)}=(\exp(\kappa\sqrt{\frac{T}{n}})+1)^{-1}$ and $1-p^{(n)}=(\exp(-\kappa\sqrt{\frac{T}{n}})+1)^{-1}$, respectively. Namely, we consider CRR markets where stock prices $S_m=S_{m/n}^{(n)}(z)$, $m=0,1,2,\dots$, satisfy (1.2) with $\rho_k=\rho_k^n$ given by (1.8) and, in addition, in place of the interest rate $r$ in the first formula in (1.4), we take the sequence of interest rates $r_n=\exp(rT/n)-1$, where $r$ is the interest rate of the BS market appearing in the second formula of (1.4) and in (1.6). We consider $S^{(n)}(z)=S^{(n)}(z,\omega)$ as a random function on $[0,T]$, so that $S^{(n)}(z,\omega)\in M[0,T]$ takes the value $S_t^{(n)}(z)=S_t^{(n)}(z,\omega)$ at $t\in[0,T]$. For $k=0,1,2,\dots,n$, put

$$
Y_k=Y_k^{(n)}(z)=F_{kT/n}(S^{(n)}(z))\quad\text{and}\quad X_k=X_k^{(n)}(z)=G_{kT/n}(S^{(n)}(z)).
\tag{2.6}
$$

Then for each $n$, the fair price $V=V^{(n)}(z)$ of the game option in the corresponding CRR market with an initial value $z>0$ of the stock is given by (1.5). Set also

$$
\begin{aligned}
&\hat S_t^{(n)}(z)=z\exp\left(\sum_{k=1}^{[nt/T]}\left(\left(r-\frac{\kappa^2}{2}\right)\frac{T}{n}+\kappa\left(\frac{T}{n}\right)^{1/2}\hat\xi_k\right)\right),\qquad t\ge T/n\quad\text{and}\\
&S_t^{(n)}=S_0^{(n)}=z>0,\qquad t\in[0,T/n),
\end{aligned}
\tag{2.7}
$$

where $\hat\xi_1,\hat\xi_2,\dots$ are i.i.d. random variables such that $\hat\xi_1=1$ or $\hat\xi_1=-1$ with the same probability 1/2 and the corresponding product measure we denote by $P^{\hat\xi}$.

Set

$$
\begin{aligned}
R_z^B(s,t)&=F_t(S^B(z))\mathbb{I}_{s\ge t}+G_s(S^B(z))\mathbb{I}_{s<t},\\
Q_z^B(s,t)&=e^{-rs\wedge t}R_z^B(s,t),
\end{aligned}
\tag{2.8}
$$

$$
\begin{aligned}
R_z^{(n)}(s,t)&=F_t(S^{(n)}(z))\mathbb{I}_{s\ge t}+G_s(S^{(n)}(z))\mathbb{I}_{s<t},\\
Q_z^{(n)}(s,t)&=e^{-rs\wedge t}R_z^{(n)}(s,t),
\end{aligned}
\tag{2.9}
$$

and let $\hat{R}_z^{(n)}(s,t), \hat{Q}_z^{(n)}(s,t)$ be defined by (2.9) with $\hat{S}^{(n)}(z)$ in place of $S^{(n)}(z)$. Denote by $\mathcal{T}_{0T}^B$, $\mathcal{T}_{0n}^\xi$ and $\mathcal{T}_{0n}^{\hat{\xi}}$ the sets of stopping times with respect to the Brownian filtration $\mathcal{F}_t^B$, $t \geq 0$, with values in $[0,T]$ and with respect to the filtrations $\mathcal{F}_k^\xi = \sigma\{\xi_1, \ldots, \xi_k\}$ and $\mathcal{F}_k^{\hat{\xi}} = \sigma\{\hat{\xi}_1, \ldots, \hat{\xi}_k\}$, $k = 0, 1, 2, \ldots$, respectively, with values in $\{0, 1, \ldots, n\}$. Set

$$V(z) = \inf_{\sigma \in \mathcal{T}_{0T}^B} \sup_{\tau \in \mathcal{T}_{0T}^B} E^B Q_z^B(\sigma, \tau), \tag{2.10}$$

$$V^{(n)}(z) = \min_{\zeta \in \mathcal{T}_{0n}^\xi} \max_{\eta \in \mathcal{T}_{0n}^\xi} E_n^\xi Q_z^{(n)}\left(\frac{\zeta T}{n}, \frac{\eta T}{n}\right) \tag{2.11}$$

and

$$\hat{V}^{(n)}(z) = \min_{\zeta \in \mathcal{T}_{0n}^\xi} \max_{\eta \in \mathcal{T}_{0n}^\xi} E^{\hat{\xi}} \hat{Q}_z^{(n)}\left(\frac{\zeta T}{n}, \frac{\eta T}{n}\right), \tag{2.12}$$

where $E^B$, $E_n^\xi$ and $E^{\hat{\xi}}$ are the expectations with respect to the probability measures $P^B$, $P_n^\xi$ and $P^{\hat{\xi}}$, respectively, and we observe that $\mathcal{T}_{0n}^\xi$ and $\mathcal{T}_{0n}^{\hat{\xi}}$ are finite sets so that we can use min and max in (2.11) and (2.12).

Recall that we choose $P^B$ to be the martingale measure for the BS market and observe that $P_n^\xi$ is the martingale measure for the corresponding CRR market since a direct computation shows that $E_n^\xi \rho_k^{(n)} = r_n$. Thus, (2.10) and (2.11) give fair prices of the game options in the corresponding markets. On the other hand, $P^{\hat{\xi}}$ is not a martingale measure, in general, and so (2.12) gives the price of the Dynkin game, but not the fair price of the corresponding game option. We note also that all our formulas involving the expectations $E^B$, in particular, (2.10) giving the fair price $V$ of a game option, do not depend on a particular choice of a continuous in time version of the Brownian motion since all of them induce the same probability measure on the space of continuous sample paths (see, e.g., Chapter 2 in [29]) which already determines all expressions with the expectations $E^B$ appearing in this paper.

The following result provides an estimate for the error term in approximation of the fair price of a game option in the BS market by fair prices of the sequence of game options and prices of Dynkin's games defined above.

THEOREM 2.1. *Suppose that $V(z)$ and $V^{(n)}(z)$ are defined by (2.9)–(2.12) with functions $F$ and $G = F + \Delta$ satisfying (2.1) and (2.2). Then there exists a constant $C > 0$ (which is, essentially, explicitly estimated in the proof) such that*

$$\begin{aligned} &\max(|V(z) - V^{(n)}(z)|, |V(z) - \hat{V}^{(n)}(z)|) \\ &\quad \leq C(F_0(z) + \Delta_0(z) + z + 1) n^{-1/4} (\ln n)^{3/4} \end{aligned} \tag{2.13}$$

*for all* $z, n > 0$.

The estimates of Theorem 2.1 remain true with, essentially, the same proof if we define $\hat{V}^{(n)}(z)$ by (2.7) and (2.12) with $\hat{\xi}_1, \hat{\xi}_2, \dots$ being arbitrary i.i.d. bounded random variables such that $E\hat{\xi}_i = 0$ and $E\hat{\xi}_i^2 = 1$. We can choose more general i.i.d. random variables $\xi_1, \xi_2, \dots$ appearing in the definition of $V^{(n)}$, as well, but these generalizations do not seem to have a financial mathematics motivation since we want to approximate game options in the BS market by simplest possible models which are, of course, game options in the CRR market.

Among main examples of options with path-dependent payoff, we have in mind integral options where

$$F_t(\upsilon) = \left(\int_0^t f_u(\upsilon_u)\, du - L\right)^+ \qquad \text{(call option case)}$$

or

$$F_t(\upsilon) = \left(L - \int_0^t f_u(\upsilon_u)\, du\right)^+ \qquad \text{(put option case)},$$

where, as usual, $a^+ = \max(a, 0)$. The penalty functional may also have here the integral form

$$\Delta_t(\upsilon) = \int_0^t \delta_u(\upsilon_u)\, du.$$

In order to satisfy conditions (2.1) and (2.2), we can assume that, for some $K > 0$ and all $x, y, u$,

$$|f_u(x) - f_u(y)| + |\delta_u(x) - \delta_u(y)| \le K|x - y|$$

and

$$|f_u(x)| + |\delta_u(x)| \le K|x|.$$

Observe also that the Asian type (averaged integral) payoffs of the form

$$F_t(\upsilon) = \left(\frac{1}{t}\int_0^t f_u(\upsilon_u)\, du - L\right)^+ \quad \text{or} \quad = \left(L - \frac{1}{t}\int_0^t f_u(\upsilon_u)\, du\right)^+$$

do not satisfy condition (2.2) if arbitrarily small exercise times are allowed, though the latter seems to have only some theoretical interest, as it hardly happens in reality. Still, also in this case, the binomial approximation errors can be estimated in a similar way considering separately estimates for small stopping times and for stopping times bounded away from zero. Namely, define $V_\varepsilon(z)$ and $V_\varepsilon^{(n)}(z)$ for $\varepsilon \ge 0$ by (2.10) and (2.11), where $Q_z^{(B)}(\sigma, \tau)$ and $Q_z^{(n)}(\frac{\zeta T}{n}, \frac{\eta T}{n})$ are replaced by $Q_z^{(B)}(\sigma \vee \varepsilon, \tau \vee \varepsilon)$ and $Q_z^{(n)}(\frac{\zeta T}{n} \vee \varepsilon, \frac{\eta T}{n} \vee \varepsilon)$, respectively. Assuming that $f_u$ and $\delta_u$ are Lipschitz continuous also in $u$ (at least for $u$ close to 0) in the form $|f_s(x) - f_u(x)| + |\delta_s(x) - \delta_u(x)| \leq K(x+1)|s-u|$ for some $K > 0$ and all $s, u, x \geq 0$, we obtain that if $\upsilon_0 = z$ and $F_0(\upsilon) = (f_0(z) - L)^+$ or $= (L - f_0(z))^+$, then

$$|F_s(\upsilon) - F_0(z)| \leq Ks\Big(1 + \sup_{0\leq u\leq s} |\upsilon_u|\Big) + K \sup_{0\leq u\leq s} |\upsilon_u - z|.$$

Using some of the estimates of Section 4, it is not difficult to see from here that $|V(z) - V_\varepsilon(z)|$ and $|V^{(n)}(z) - V_\varepsilon^{(n)}(z)|$ do not exceed $C(1+z)\sqrt{\varepsilon}$ for all small $\varepsilon$ and some constant $C$. On the other hand, applying the same estimates as in the proof of Theorem 2.1, we derive that, for some constant $C > 0$ and all $n, \varepsilon > 0$,

$$|V_\varepsilon(z) - V_\varepsilon^{(n)}(z)| \leq C(1+z)\varepsilon^{-1} n^{-1/4} (\ln n)^{3/4}.$$

Choosing $\varepsilon = n^{-1/6}\sqrt{\ln n}$, we obtain that, under the above conditions in the case of Asian options, $|V(z) - V^{(n)}(z)|$ can be estimated by $3C(1+z)n^{-1/12}(\ln n)^{1/4}$.

Another important example of path-dependent payoffs are, so-called, Russian options where, for instance,

$$F_t(\upsilon) = \max\Big(m, \sup_{u\in[0,t]} \upsilon_u\Big) \quad \text{and} \quad \Delta_t(\upsilon) = \delta\upsilon_t.$$

Such payoffs satisfy the conditions of Theorem 2.1. Indeed, (2.1) is clear in this case and (2.2) follows since, for $t \geq s$,

$$\begin{aligned}\max\Big(m, \sup_{u\in[0,t]} \upsilon_u\Big) - \max\Big(m, \sup_{u\in[0,s]} \upsilon_u\Big) &\leq \sup_{u\in[0,t]} \upsilon_u - \sup_{u\in[0,s]} \upsilon_u \\ &\leq \sup_{u\in[s,t]} \upsilon_u - \upsilon_s \\ &\leq \sup_{u\in[s,t]} |\upsilon_u - \upsilon_s|.\end{aligned}$$

In fact, the estimates (2.13) can be improved a bit for the Russian options case dropping the logarithmic term there (see Remark 3.7). Of course, conditions (2.1) and (2.2) are always satisfied in the case of standard options with payoffs depending only on the current stock price.

Observe that many other path-dependent and look back options have payoffs which can be represented via functions $F$ and $\Delta$ satisfying (2.1) and (2.2). Barrier options have discontinuous payoffs which cannot be represented via functions satisfying (2.1) and (2.2), but a small modification of our approach goes through in this case as well. This modification is based on the observation that troubles with the approximation occur here when the supremum of stock prices (or a similar quantity) belongs to a small neighborhood of the barrier value, but the probability of this event is small since this supremum is usually a random variable having a bounded probability density function.

In order to compare $V(z)$ and $V^{(n)}(z)$ in the case of path dependent payoffs, we have to consider both BS and CRR markets on one probability space in an appropriate way and the main tool in achieving this goal will be here the Skorokhod type embedding (see, e.g., [4], Section 37). In fact, for the binomial i.i.d. random variables $\xi_1, \xi_2, \ldots$ and $\hat{\xi}_1, \hat{\xi}_2, \ldots$ appearing in the setup of the CRR market models above, the embedding is explicit and no general theorems are required, but if we want to extend the result for other sequences of i.i.d. random variables, we have to rely upon the general result. Namely, define recursively

$$\theta_0^{(n)} = 0, \qquad \theta_{k+1}^{(n)} = \inf\left\{t > \theta_k^{(n)} : |B_t^* - B^*_{\theta_k^{(n)}}| = \sqrt{\frac{T}{n}}\right\}$$

where, recall, $B_t^* = -\frac{\kappa}{2}t + B_t$ and

$$\hat{\theta}_0^{(n)} = 0, \qquad \hat{\theta}_{k+1}^{(n)} = \inf\left\{t > \hat{\theta}_k^{(n)} : |B_t - B_{\hat{\theta}_k^{(n)}}| = \sqrt{\frac{T}{n}}\right\}.$$

The standard strong Markov property based arguments (cf. [4], Section 37) show that $\theta_k^{(n)} - \theta_{k-1}^{(n)}, k = 1, 2, \ldots,$ and $\hat{\theta}_k^{(n)} - \hat{\theta}_{k-1}^{(n)}$, $k = 1, 2, \ldots,$ are i.i.d. sequences of random variables such that $(\theta_{k+1}^{(n)} - \theta_k^{(n)}, B^*_{\theta_{k+1}^{(n)}} - B^*_{\theta_k^{(n)}})$ and $(\hat{\theta}_{k+1}^{(n)} - \hat{\theta}_k^{(n)}, B_{\hat{\theta}_{k+1}^{(n)}} - B_{\hat{\theta}_k^{(n)}})$ are independent of $\mathcal{F}^B_{\theta_k^{(n)}}$ (where, recall, $\mathcal{F}_t^B = \sigma\{B_s, s \le t\}$). This is standard for the Brownian motion $B_t$ and the stopping times $\hat{\theta}_k^{(n)}$ (see [4], Section 37), but can be proved by exactly the same method for the Brownian motion with a constant drift $B_t^*$ and the stopping times $\theta_k^{(n)}$ as well.

Another way to justify this independency assertion is to apply the Girsanov theorem (see [11, 13] and [28]) which is useful to have in mind in our situation anyway. Namely, for any Brownian stopping time $\tau$ satisfying $E^B \exp(\kappa^2\tau/8) < \infty$, the process $B_t^*$, $t \ge 0$, becomes a standard Brownian motion on the interval $[0, \tau]$ with respect to the probability measure $P_\tau^*$ determined by

$$dP_\tau^* = q(\tau, B_\tau)\, dP^B \quad \text{or} \quad dP^B = q^*(\tau, B_\tau^*)\, dP_\tau^*, \tag{2.14}$$

where

$$q(t, B_t) = \exp\left(\frac{\kappa}{2}B_t - \frac{1}{8}\kappa^2 t\right) = \exp\left(\frac{\kappa}{2}B_t^* + \frac{1}{8}\kappa^2 t\right) = (q^*(t, B_t^*))^{-1}.$$

Observe that if $\sigma > \tau$ is another Brownian stopping time, then $P_\sigma^* = P_\tau^*$ on the $\sigma$-algebra $\mathcal{F}_\tau^B$ so the above probability measures enable us to consider $B^*_t$ as a Brownian motion on any time interval. Hence, under the measure $P^*$, the process $B^*_t$, $t\ge 0$, together with the stopping times $\theta^{(n)}_k$, $k\ge 0$, becomes the Brownian motion $B_t$, $t\ge 0$, together with the stopping times $\hat\theta^{(n)}_k, k\ge 0$. Hence, by the strong Markov property of the Brownian motion, any two events $A\in\mathcal{F}^B_{\theta^{(n)}_k}$ and $\tilde A\in\sigma\{B^*_{t+\theta^{(n)}_k}-B^*_{\theta^{(n)}_k}, t\in[0,\theta^{(n)}_{k+1}-\theta^{(n)}_k]\}$, are independent under $P^*_{\theta^{(n)}_{k+1}}$. But then $A$ and $\tilde A$ will remain independent under the original probability $P^B$ in view of the factorization property of the density $q^*$,

$$q^*(\theta^{(n)}_{k+l},B^*_{\theta^{(n)}_{k+l}})=q^*(\theta^{(n)}_k,B^*_{\theta^{(n)}_k})q^*(\theta^{(n)}_{k+l}-\theta^{(n)}_k,B^*_{\theta^{(n)}_{k+l}}-B^*_{\theta^{(n)}_k})$$

with the first factor measurable with respect to $\mathcal{F}^B_{\theta^{(n)}_k}$ and the second factor independent of it and having the same distribution as $q^*(\theta^{(n)}_l,B^*_{\theta^{(n)}_l})$. It follows that $(\theta^{(n)}_{k+1}-\theta^{(n)}_k,B^*_{\theta^{(n)}_{k+1}}-B^*_{\theta^{(n)}_k})$ is independent of $\mathcal{F}^B_{\theta^{(n)}_k}$ under the original probability $P^B$ as well.

It turns out (see [30] and the beginning of Section 4) that $B^*_{\theta^{(n)}_1}$ has the same distribution as $\sqrt{\frac{T}{n}}\xi_1$. Clearly, $B_{\hat\theta^{(n)}_1}=\sqrt{\frac{T}{n}}$ or $=-\sqrt{\frac{T}{n}}$ with the same probability 1/2, and so $B_{\hat\theta^{(n)}_1}$ has the same distribution as $\sqrt{\frac{T}{n}}\hat\xi_1$. Set

$$\Xi^{(n)}_k=\left(\frac{T}{n}\right)^{1/2}\sum_{j=1}^k\xi_j\quad\text{and}\quad\hat\Xi^{(n)}_k=\left(\frac{T}{n}\right)^{1/2}\sum_{j=1}^k\hat\xi_j,\tag{2.15}$$

then $\Xi^{(n)}_k$ and $\hat\Xi^{(n)}_k$ have the same distribution as $B^*_{\theta^{(n)}_k}$ and $B_{\hat\theta^{(n)}_k}$, respectively.

Theorem 2.1 provides an approximation of the fair price of game options in the BS market by means of fair prices of game options in the CRR market which becomes especially useful if we can provide also a simple description of rational (or $\delta$-rational) exercise times of these options in the BS market via exercise times of their CRR market approximations which are, by the definition, optimal (or $\delta$-optimal) stopping times for the Dynkin games whose price are given by (2.11) and (2.12), respectively. For each $k=1,2,\dots,$ introduce the finite $\sigma$-algebra $\mathcal{G}^{B,n}_k=\sigma\{B^*_{\theta^{(n)}_1},B^*_{\theta^{(n)}_2}-B^*_{\theta^{(n)}_1},\dots,B^*_{\theta^{(n)}_k}-B^*_{\theta^{(n)}_{k-1}}\}$, which is, clearly, isomorphic to $\mathcal{F}^\xi_k=\sigma\{\xi_i,i\le k\}$ considered before since each element of $\mathcal{G}^{B,n}_k$ and of $\mathcal{F}^\xi_k$ is an event of the form

$$A^{B,n}_{\iota^{(k)}}=\left\{B^*_{\theta^{(n)}_j}-B^*_{\theta^{(n)}_{j-1}}=\iota_j\sqrt{\frac{T}{n}},j=1,\dots,k\right\}$$

and

$$A^{\xi}_{\iota^{(k)}} = \{\xi_j = \iota_j, j = 1, \dots, k\},$$

respectively, where $\iota^{(k)} = (\iota_1, \dots, \iota_k) \in \{-1,1\}^k$, $\theta_0^{(n)} = 0$ and $B_0 = 0$. Let $\mathcal{S}^{B,n}$ be the set of stopping times with respect to the filtration $\mathcal{G}_k^{B,n}, k = 0, 1, 2, \dots,$ where $\mathcal{G}_0^{B,n} = \{\varnothing, \Omega_B\}$ is the trivial $\sigma$-algebra and $\Omega_B$ is the sample space of the Brownian motion. The subset of these stopping times with values in $\{0, 1, \dots, n\}$ will be denoted by $\mathcal{S}_{0,n}^{B,n}$. For each $\iota^{(n)} = (\iota_1, \dots, \iota_n) \in \{-1,1\}^n$ and $k < n$, we set $\iota^{(k)} = (\iota_1, \dots, \iota_k) \in \{-1,1\}^k$. Denote by $\mathcal{J}_{0,n}$ the set of functions $\nu : \{-1,1\}^n \to \{0, 1, \dots, n\}$ such that if $\nu(\iota^{(n)}) = k \le n$ and $\tilde{\iota}^{(k)} = \iota^{(k)}$ for some $\tilde{\iota}^{(n)} \in \{-1,1\}^n$, then $\nu(\tilde{\iota}^{(n)}) = k$ as well. Define the functions $\lambda_\xi^{(n)} : \Omega_\xi \to \{-1,1\}^n$ and $\lambda_B^{(n)} : \Omega_B \to \{-1,1\}^n$ by $\lambda_\xi^{(n)}(\omega) = (\xi_1(\omega), \dots, \xi_n(\omega))$ and

$$\lambda_B^{(n)}(\omega) = \sqrt{\frac{n}{T}}\big(B^*_{\theta_1^{(n)}(\omega)}(\omega), B^*_{\theta_2^{(n)}(\omega)}(\omega) - B^*_{\theta_1^{(n)}(\omega)}(\omega), \dots, \\ B^*_{\theta_n^{(n)}(\omega)}(\omega) - B^*_{\theta_{n-1}^{(n)}(\omega)}(\omega)\big),$$

where $\Omega_\xi$ and $\Omega_B$ are sample spaces on which the sequence $\xi_1, \xi_2, \dots$ and the Brownian motion $B_t$ are defined, respectively. It is clear that any $\zeta \in \mathcal{T}_{0n}^{\xi}$ and $\eta \in \mathcal{S}_{0,n}^{B,n}$ can be represented uniquely in the form $\zeta = \mu \circ \lambda_\xi^{(n)}$ and $\eta = \nu \circ \lambda_B^{(n)}$ for some $\mu, \nu \in \mathcal{J}_{0,n}$. Similarly, we introduce $\hat{\mathcal{G}}_k^{B,n} = \sigma\{B_{\hat{\theta}_1^{(n)}}, B_{\hat{\theta}_2^{(n)}} - B_{\hat{\theta}_1^{(n)}}, \dots, B_{\hat{\theta}_k^{(n)}} - B_{\hat{\theta}_{k-1}^{(n)}}\}$, which is isomorphic to $\hat{\mathcal{F}}_k^{\xi} = \sigma\{\hat{\xi}_i, i \le k\}$, $\hat{\lambda}_\xi^{(n)}(\omega) = (\hat{\xi}_1(\omega), \dots, \hat{\xi}_n(\omega))$ and

$$\hat{\lambda}_B^{(n)}(\omega) = \sqrt{\frac{n}{T}}\big(B_{\hat{\theta}_1^{(n)}(\omega)}(\omega), B_{\hat{\theta}_2^{(n)}(\omega)}(\omega) - B_{\hat{\theta}_1^{(n)}(\omega)}(\omega), \dots, \\ B_{\hat{\theta}_n^{(n)}(\omega)}(\omega) - B_{\hat{\theta}_{n-1}^{(n)}(\omega)}(\omega)\big).$$

The following result will be proved in Section 5.

THEOREM 2.2. *There exists a constant $C > 0$ (which is, essentially, estimated explicitly in the proof) such that if $\zeta_n^* = \mu_n^* \circ \lambda_\xi^{(n)}$ and $\eta_n^* = \nu_n^* \circ \lambda_\xi^{(n)}$, $\mu_n^*, \nu_n^* \in \mathcal{J}_{0n}$, are rational exercise times for the game option in the CRR market defined by (2.5), that is,*

$$V^{(n)}(z) = \min_{\zeta \in \mathcal{T}_{0n}^{\xi}} E^{\xi} Q_z^{(n)}\left(\zeta\frac{T}{n}, \eta_n^*\frac{T}{n}\right) = \max_{\eta \in \mathcal{T}_{0n}^{\xi}} E^{\xi} Q_z^{(n)}\left(\zeta_n^*\frac{T}{n}, \eta\frac{T}{n}\right), \tag{2.16}$$

*then* $\varphi_n^*=\theta^{(n)}_{\mu_n^*\circ\lambda_B^{(n)}}$ *and* $\psi_n^*=\theta^{(n)}_{\nu_n^*\circ\lambda_B^{(n)}}$ *are* $\delta_n(z)$*-rational exercise times for the game option in the BS market defined by* (2.3) *and* (2.4)*, that is,*

$$\sup_{\tau\in\mathcal{T}_{0T}^B} E^B Q_z^B(\varphi_n^*,\tau)-\delta_n(z)\le V(z)\le \inf_{\sigma\in\mathcal{T}_{0T}^B} E^B Q_z^B(\sigma,\psi_n^*)+\delta_n(z), \tag{2.17}$$

*where* $\delta_n(z)=C(F_0(z)+\Delta_0(z)+z+1)n^{-1/4}(\ln n)^{3/4}$. *The assertions remain true if we replace above* $\lambda_\xi^{(n)}$, $\lambda_B^{(n)}$, $Q_z^{(n)}$ *and* $V^{(n)}(z)$ *by* $\hat\lambda_\xi^{(n)}$, $\hat\lambda_B^{(n)}$, $\hat Q_z^{(n)}$ *and* $\hat V^{(n)}(z)$*, respectively.*

It is well known (see, e.g., [26]) that when payoffs depend only on the current stock price (a Markov case), $\delta$-optimal stopping times of Dynkin's games can be obtained as first arrival times to domains where the payoff is $\delta$-close to the value of the game (as a function of the initial stock price). For path dependent payoffs, the situation is more complicated and, in general, in order to construct $\delta$-optimal stopping times, we have to know the stochastic process of values of the games starting at each time $t\in[0,T]$ conditioned to the information up to $t$. It is not clear what kind of approximation of this process can provide some information about $\delta$-rational exercise times and the convenient alternative method of their construction exhibited in Theorem 2.2 seems to be important both for the theory and applications. Moreover, this construction is effective and can be employed in practice since $\mu_n^*$ and $\nu_n^*$ are functions on sequences of 1's and $-1$'s which can be computed (and stored in a computer) using the recursive formulas (1.7) even before the stock evolution begins. In order to compute $\lambda_B^{(n)}$, we have to watch the discounted stock price $\check S_t^B(z)=e^{-rt}S_t^B(z)$ evolution of a real stock at moments $\theta_k^{(n)}$ which are obtained recursively by $\theta_0^{(n)}=0$ and

$$\theta_{k+1}^{(n)}=\inf\{t>\theta_k^{(n)}:\check S_t^B(z)=e^{\pm\kappa(T/n)^{1/2}}\check S^B_{\theta_k^{(n)}}(z)\} \tag{2.18}$$

and to construct the $\{1,-1\}$ sequence $\lambda_B^{(n)}(\omega)$ by writing 1 or $-1$ on $k$th place depending on whether $\check S^B_{\theta_k^{(n)}}(z)=e^{\kappa(T/n)^{1/2}}\check S^B_{\theta_{k-1}^{(n)}}(z)$ or $\check S^B_{\theta_k^{(n)}}(z)=e^{-\kappa(T/n)^{1/2}}\check S^B_{\theta_{k-1}^{(n)}}(z)$, respectively. Observe also that Theorem 2.2 can be extended to more general sequences of random variables $\xi_1,\xi_2,\dots$ and $\hat\xi_1,\hat\xi_2,\dots$, but this does not seem to have much of an interest for applications.

Recall (see [28]) that a sequence $\pi=(\pi_1,\dots,\pi_n)$ of pairs $\pi_k=(\beta_k,\gamma_k)$ of $\mathcal{F}_{k-1}^\xi$-measurable random variables $\beta_k$, $\gamma_k$, $k=1,\dots,n$, is called a self-financing portfolio strategy in the CRR market determined by (1.2), (1.4), (1.8) and (2.5) if the price of the portfolio at time $k$ is given by the formula

$$Z_k^{\pi,n}=\beta_k b_k+\gamma_k S_{kT/n}^{(n)}(z)=\beta_{k+1}b_k+\gamma_{k+1}S_{kT/n}^{(n)}(z) \tag{2.19}$$

and the latter equality means that all changes in the portfolio value are due to capital gains and losses but not to withdrawal or infusion of funds. A pair $(\zeta,\pi)$ of a stopping time $\zeta\in\mathcal{T}_{0n}^{\xi}$ and a self-financing portfolio strategy $\pi$ is called a hedge for (against) the game option with the payoff $R_z^{(n)}$ given by (2.9) if (see [15])

$$Z_{\zeta\wedge k}^{\pi;n}\geq R_z^{(n)}\left(\frac{\zeta T}{n},\frac{kT}{n}\right)\qquad\forall\,k=0,1,\dots,n. \tag{2.20}$$

It follows from [15] that, for any $\zeta\in\mathcal{T}_{0n}^{\xi}$, there exists a self-financing portfolio strategy $\pi^{\zeta}$ so that $(\zeta,\pi^{\zeta})$ is a hedge. In particular, if we take the rational exercise time $\zeta=\zeta_n^*$ of the writer, then such $\pi^{\zeta}$ exists with the initial portfolio capital $V^{(n)}(z)$. The construction of $\pi^{\zeta}$ goes directly via the Doob decomposition of supermartingales and a martingale representation lemma (see [15] and [28]), both being explicit in the CRR market case. In the continuous time BS market we cannot write the corresponding portfolio strategies in an explicit way, and so some approximations are necessary though, surprisingly, this problem has not been treated before in the literature.

THEOREM 2.3. *Let $\zeta\in\mathcal{T}_{0n}^{\xi}$, $\pi=\pi^{\zeta}$ and (2.19) together with (2.20) hold true with $\mathcal{F}_k^{\xi}$-measurable $\beta_k=\beta_k^{\zeta}$ and $\gamma_k=\gamma_k^{\xi}$, so that $(\zeta,\pi^{\zeta})$ is a hedge. Then $\beta_k^{\zeta}=f_k\circ\lambda_{\xi}^{(k-1)}$, $\gamma_k^{\zeta}=g_k\circ\lambda_{\xi}^{(k-1)}$ and $\zeta=\mu\circ\lambda_{\xi}^{(n)}$ for some uniquely defined functions $f_k$, $g_k$ on $\{-1,1\}^{k-1}$ and some $\mu\in\mathcal{J}_{0n}$. Let $\varphi=\mu\circ\lambda_B^{(n)}$ and set $\beta_t^{\varphi}=f_k\circ\lambda_B^{(k-1)}$ and $\gamma_t^{\varphi}=g_k\circ\lambda_B^{(k-1)}$ whenever $t\in(\theta_{k-1}^{(n)},\theta_k^{(n)}]$. Then*

$$Z_t^B=\beta_t^{\varphi}b_t+\gamma_t^{\varphi}S_t^B(z) \tag{2.21}$$

*is a self-financing portfolio in the BS market and there exists a constant $C>0$ such that*

$$\begin{aligned} &E^B\sup_{0\leq t\leq T}(R_z^B(\theta_{\varphi}^{(n)},t)-Z^B_{\theta_{\varphi}^{(n)}\wedge t})^+\\ &\qquad\leq C(F_0(z)+\Delta_0(z)+z+1)n^{-1/4}(\ln n)^{3/4},\end{aligned} \tag{2.22}$$

*where $a^+=\max(a,0)$. In particular, there exists a self-financing portfolio of this form satisfying (2.22) with the initial value $V^{(n)}(z)$ [which according to (2.13) is close to the fair price $V(z)$ of the game option] taking $\varphi^*=\mu^*\circ\lambda_B^{(n)}$ if $\zeta^*=\mu^*\circ\lambda_{\xi}^{(n)}$ is the rational exercise time and $\pi=\pi^{\zeta^*}$ is the corresponding optimal self-financing hedging portfolio strategy for the CRR market.*

Inequality (2.22) estimates the expectation of the maximal shortfall (risk) of certain (nearly hedging) portfolio strategy which can be constructed effectively in applications since the functions $f_l,g_l$ and $\mu$ are determined by a self-financing hedging strategy in the CRR market which can be computed directly and stored in a computer even before the real stock evolution begins or in case of computer memory limitations, we can compute these functions each time when needed using corresponding algorithms for the CRR market. The functions $\lambda_B^{(n)}$ or, in other words, the sequences from $\{-1,1\}^n$ which should be plugged into the functions $f_l$, $g_l$ and $\mu$ should be obtained in practice by watching the evolution of the discounted stock price $e^{-rt}S_t^B$ at moments $\theta_k^{(n)}$ as described after Theorem 2.2.

**3. Auxiliary lemmas.** In addition to the set $\mathcal{S}^{B,n}$ of stopping times with respect to the filtration $\{\mathcal{G}_k^{B,n}\}_{k=0,1,2,\dots}$ introduced before Theorem 2.2, consider also the set $\mathcal{T}^{B,n}$ of stopping times with respect to the filtration $\{\mathcal{F}^B_{\theta_k^{(n)}}\}_{k=0,1,2,\dots}$ with values in $\{0,1,2,\dots\}$ and the subset of such stopping times with values in $\{0,1,\dots,n\}$ will be denoted by $\mathcal{T}_{0,n}^{B,n}$. Clearly, $\mathcal{S}^{B,n}\subset\mathcal{T}^{B,n}$. Set

$$S_t^{B,n}(z)=z\exp\left(\sum_{k=1}^{[nt/T]}\left(\frac{rT}{n}+\kappa(B^*_{\theta_k^{(n)}}-B^*_{\theta_{k-1}^{(n)}})\right)\right)\qquad\text{if } t\in[T/n,T],$$

$$S_t^{B,n}=S_T^{B,n}\qquad\text{if } t>T\quad\text{and}\tag{3.1}$$

$$S_t^{B,n}=S_0^{B,n}=z>0\qquad\text{if } t\in[0,T/n)$$

and let $\hat S_t^{B,n}$ be the corresponding expression if we replace in (3.1) $r$ by $r-\frac{\kappa^2}{2}$, $B^*$ by $B$, and $\theta_k^{(n)}$'s by $\hat\theta_k^{(n)}$. Denote

$$\begin{aligned}R_z^{B,n}(s,t)&=F_t(S^{B,n}(z))\mathbb{I}_{s\ge t}+G_s(S^{B,n}(z))\mathbb{I}_{s<t},\\ Q_z^{B,n}(s,t)&=e^{-rs\wedge t}R_z^{B,n}(s,t),\end{aligned}\tag{3.2}$$

$$V^{B,n}(z)=\inf_{\zeta\in\mathcal{T}_{0,n}^{B,n}}\sup_{\eta\in\mathcal{T}_{0,n}^{B,n}}E^BQ_z^{B,n}\left(\frac{\zeta T}{n},\frac{\eta T}{n}\right)\tag{3.3}$$

and

$$V_{\mathcal{S}}^{B,n}(z)=\min_{\zeta\in\mathcal{S}_{0,n}^{B,n}}\max_{\eta\in\mathcal{S}_{0,n}^{B,n}}E^BQ_z^{B,n}\left(\frac{\zeta T}{n},\frac{\eta T}{n}\right).\tag{3.4}$$

Consider also the corresponding quantities $\hat R_z^{B,n}(s,t)$, $\hat Q_z^{B,n}(s,t)$, $\hat V^{B,n}(z)$ and $\hat V_{\mathcal{S}}^{B,n}(z)$ which are obtained by taking in the above formulas $\hat S_t^{B,n}$ in place of $S_t^{B,n}$. Though $V_{\mathcal{S}}^{B,n}(z)$ and $\hat V_{\mathcal{S}}^{B,n}(z)$ are not used in the proofs, their introduction clarifies the nature of various sets of stopping times involved here.

The reason for considering the Skorokhod embedding here and the basis for our proofs of Theorems 2.1–2.3 is the following result which is a generalization of Lemma 3.1 from [24] and which enables us to consider all relevant processes on one probability space and to deal with stopping times with respect to the same Brownian filtration only.

LEMMA 3.1. *For any* $z, n > 0$,

$$V_{\mathcal{S}}^{B,n}(z) = V^{(n)}(z) = V^{B,n}(z). \tag{3.5}$$

*The same result holds true if we replace in (3.5) all $V$'s by $\hat{V}$'s.*

PROOF. First, observe that $\theta_\zeta^{(n)} \in \mathcal{T}^B$ for any $\zeta \in \mathcal{T}^{B,n}$ (see [24]), where $\mathcal{T}^B$ is the set of all almost surely (a.s.) finite stopping times for the Brownian motion $B_t, t \geq 0$. Indeed,

$$\{\theta_\zeta^{(n)} \leq t\} = \bigcup_{k=0}^{n} \{\theta_k^{(n)} \leq t\} \cap \{\zeta = k\} \tag{3.6}$$

and since $\{\zeta = k\} \in \mathcal{F}^B_{\theta_k^{(n)}}$ and $\{\theta_k^{(n)} \leq t\} \in \mathcal{F}_t^B$, we conclude that the event in the right-hand side of (3.6) belongs to $\mathcal{F}_t^B$, and so $\theta_\zeta^{(n)}$ is a stopping time.

Next, as we mentioned before the statement of Theorem 2.2, $\mathcal{T}_{0n}^{\xi} = \{\mu \circ \lambda_\xi^{(n)} : \mu \in \mathcal{J}_{0,n}\}$ and $\mathcal{S}_{0,n}^{B,n} = \{\nu \circ \lambda_B^{(n)} : \nu \in \mathcal{J}_{0,n}\}$. It is clear that, for any $\mu, \nu \in \mathcal{J}_{0,n}$,

$$Q_z^{(n)}\left(\frac{T}{n}\mu \circ \lambda_\xi^{(n)}, \frac{T}{n}\nu \circ \lambda_\xi^{(n)}\right) \quad \text{and} \quad Q_z^{B,n}\left(\frac{T}{n}\mu \circ \lambda_B^{(n)}, \frac{T}{n}\nu \circ \lambda_B^{(n)}\right)$$

have the same distributions, and so the first equality in (3.5) follows.

In order to prove the second equality in (3.5), we employ the dynamical programming relations (1.7) for $V_{k,n}^{(n)} = V_{k,n}^{(n)}(z)$ and for $V_{k,n}^{B,n} = V_{k,n}^{B,n}(z)$, $k = 0, 1, \ldots, n$, which in our case have the form

$$\begin{aligned} &V^{(n)} = V_{0,n}^{(n)}, \qquad V_{n,n}^{(n)} = e^{-rT} F_T(S^{(n)}) \quad \text{and} \\ &V_{k,n}^{(n)} = \min(e^{-rkT/n} X_k^{(n)}(z), \max(e^{-rkT/n} Y_k^{(n)}(z), E(V_{k+1,n}^{(n)} | \mathcal{F}_k^{\xi}))) \end{aligned} \tag{3.7}$$

and

$$\begin{aligned} &V^{B,n} = V_{0,n}^{B,n}, \qquad V_{n,n}^{B,n} = e^{-rT} F_T(S^{B,n}) \quad \text{and} \\ &V_{k,n}^{B,n} = \min(e^{-rkT/n} X_k^{B,n}(z), \max(e^{-rkT/n} Y_k^{B,n}(z), E(V_{k+1,n}^{B,n} | \mathcal{F}^B_{\theta_k^{(n)}}))), \end{aligned} \tag{3.8}$$

where $X_k^{(n)}(z)$ and $Y_k^{(n)}(z)$ are given by (2.6),

$$Y_k^{B,n}(z) = F_{kT/n}(S^{B,n}(z)) \quad \text{and} \quad X_k^{B,n}(z) = G_{kT/n}(S^{B,n}(z)). \tag{3.9}$$

For any numbers $x_1, x_2, \dots, x_n$, set

$$x_t^{(n)} = x_t^{(n)}(z) = z \exp\left(\sum_{k=1}^{[nt/T]} \left(\frac{rT}{n} + \kappa x_k\right)\right) \qquad \text{if } t \ge T/n$$

and $x_t^{(n)} = x_t^{(n)}(z) = z$ if $t \in [0, T/n)$. In view of (2.1), we can write

$$\begin{aligned} F_{kT/n}(x^{(n)}(z)) &= q_k(z, x_1, \dots, x_k) \quad \text{and} \\ \Delta_{kT/n}(x^{(n)}(z)) &= r_k(z, x_1, \dots, x_k) \end{aligned} \tag{3.10}$$

for some continuous functions $q_k$ and $r_k$ depending only on $z, x_1, \dots, x_k$. Next, we show by the backward induction that there exist measurable functions $\Phi_k(z, x_1, \dots, x_k)$, $k = 1, 2, \dots, n$, and $\Phi_0(z)$ such that

$$V_{kn}^{(n)}(z) = \Phi_k\left(z, \left(\frac{T}{n}\right)^{1/2} \xi_1, \dots, \left(\frac{T}{n}\right)^{1/2} \xi_k\right), \qquad V_{0n}^{(n)}(z) = \Phi_0(z) \tag{3.11}$$

and

$$\begin{aligned} V_{kn}^{B,n}(z) &= \Phi_k(z, B^*_{\theta_1^{(n)}}, B^*_{\theta_2^{(n)}} - B^*_{\theta_1^{(n)}}, \dots, B^*_{\theta_k^{(n)}} - B^*_{\theta_{k-1}^{(n)}}), \\ V_{0n}^{B,n}(z) &= \Phi_0(z). \end{aligned} \tag{3.12}$$

Indeed, for $k = n$, set $\Phi_n(z, x_1, \dots, x_n) = e^{-rT} q_n(z, x_1, \dots, x_n)$. Suppose that the assertion holds true for $k \ge l + 1$, that is, that for such $k$'s, we found functions $\Phi_k$ satisfying (3.11) and (3.12). Now, set

$$\begin{aligned} \Phi_l(z, x_1, \dots, x_l) = \min(&e^{-rlT/n}(q_l(z, x_1, \dots, x_l) + r_l(z, x_1, \dots, x_l)), \\ &\max(e^{-rlT/n} q_l(z, x_1, \dots, x_l), h_l(z, x_1, \dots, x_l))), \end{aligned} \tag{3.13}$$

where

$$\begin{aligned} h_l(z, x_1, \dots, x_l) &= E^{\xi} \Phi_{l+1}(z, x_1, \dots, x_l, (T/n)^{1/2} \xi_{l+1}) \\ &= E^B \Phi_{l+1}(z, x_1, \dots, x_l, B^*_{\theta_{l+1}^{(n)}} - B^*_{\theta_l^{(n)}}). \end{aligned}$$

Then (3.7) and (3.8) will be satisfied for $k = l$ with $V_{kn}^{(n)}(z)$ and $V_{kn}^{B,n}(z)$ given by (3.11) and (3.12), respectively, since $\xi_{l+1}$ and $B^*_{\theta_{l+1}^{(n)}} - B^*_{\theta_l^{(n)}}$ are independent of $\mathcal{F}_l^{\xi}$ and $\mathcal{F}_l^{B,n}$, respectively, and so by the standard fact (see, e.g., Example 34.3 in [4]),

$$\begin{aligned} &E\left(\Phi_{l+1}\left(z, \left(\frac{T}{n}\right)^{1/2} \xi_1, \dots, \left(\frac{T}{n}\right)^{1/2} \xi_l, \left(\frac{T}{n}\right)^{1/2} \xi_{l+1}\right) \Big| \mathcal{F}_l^{\xi}\right) \\ &\quad = h_l\left(z, \left(\frac{T}{n}\right)^{1/2} \xi_1, \dots, \left(\frac{T}{n}\right)^{1/2} \xi_l\right) \end{aligned}$$

and

$$E^B(\Phi_{l+1}(z, B^*_{\theta_1^{(n)}}, B^*_{\theta_2^{(n)}} - B^*_{\theta_1^{(n)}}, \dots, B^*_{\theta_l^{(n)}} - B^*_{\theta_{l-1}^{(n)}}, B^*_{\theta_{l+1}^{(n)}} - B^*_{\theta_l^{(n)}})|\mathcal{F}_l^{B,n})$$
$$= h_l(z, B^*_{\theta_1^{(n)}}, B^*_{\theta_2^{(n)}} - B^*_{\theta_1^{(n)}}, \dots, B^*_{\theta_l^{(n)}} - B^*_{\theta_{l-1}^{(n)}}),$$

completing the induction. Now applying (3.11) and (3.12) with $k = 0$, we arrive at the second equality in (3.5). We obtain the assertion (3.5) for $\hat{V}$'s in place of $V$'s exactly in the same way as above. □

Next, for readers' convenience, we formulate a series of lemmas which demonstrate the plan of our proof of Theorem 2.1 leaving till the next section the actual proof of these results which rely on relatively standard stochastic analysis estimates. We will do $L^1$-estimates directly with respect to the probability $P^B$, though we could do instead $L^2$ estimates with respect to the probability $P^*_\tau$, and then pass to estimates with respect to the original measure $P^B$ using the Girsanov transformation (2.14) and the Cauchy–Schwarz inequality. This would enable us to work from the beginning with the stopping times $\hat{\theta}_1^{(n)}, \hat{\theta}_2^{(n)}, \dots$ in place of $\theta_1^{(n)}, \theta_2^{(n)}, \dots$ which is easier but, on the other hand, would need $L^2$ estimates which require few additional lines anyway. We formulate results which lead to the required estimate of $|V(z) - V^{(n)}(z)|$. The corresponding estimate of $|V(z) - \hat{V}^{(n)}(z)|$ proceeds exactly in the same way replacing $\xi_k$'s by $\hat{\xi}_k$'s and $\theta_k^{(n)}$'s by $\hat{\theta}_k^{(n)}$'s which, in fact, leads to a bit easier arguments.

First, observe that though $S_t^{B,n}(z)$ defined by (3.1) is a piecewise constant approximation of the BS stock price $S_t(z)$ given by (2.4), there is certain inconsistency there between times $kT/n$ of jumps of $S_t^{B,n}(z)$ and the values $B_{\theta_k^{(n)}} - B_{\theta_{k-1}^{(n)}}$ of jumps in the exponent which corresponds to the Brownian stopping time $\theta_k^{(n)}$. In order to pass to the correct time, we introduce

$$\begin{aligned} &S_t^{B,\theta,n}(z) = z\exp(r\theta_k^{(n)} + \kappa B^*_{\theta_k^{(n)}}) \qquad \text{if } \theta_k^{(n)} \le t < \theta_{k+1}^{(n)}, k = 0, 1, \dots, n, \\ &S_t^{B,\theta,n}(z) = S^{B,\theta,n}_{\theta_n^{(n)}}(z) \qquad \text{if } t \ge \theta_n^{(n)}. \end{aligned} \tag{3.14}$$

Set

$$\begin{aligned} &R_z^{B,\theta,n}(s,t) = F_t(S^{B,\theta,n}(z))\mathbb{I}_{s\ge t} + G_s(S^{B,\theta,n}(z))\mathbb{I}_{s<t}, \\ &Q_z^{B,\theta,n}(s,t) = e^{-rs\wedge t} R_z^{B,\theta,n}(s,t) \end{aligned} \tag{3.15}$$

and

$$V^{B,\theta,n}(z) = \inf_{\zeta \in \mathcal{T}_{0,n}^{B,n}} \sup_{\eta \in \mathcal{T}_{0,n}^{B,n}} E^B Q_z^{B,\theta,n}(\theta_\zeta^{(n)}, \theta_\eta^{(n)}). \tag{3.16}$$

In order to compare $V^{B,n}(z)$ and $V^{B,\theta,n}(z)$, we have to be able to compare $S_t^{B,n}(z)$ and $S_t^{B,\theta,n}(z)$ at the same time $t\in[0,T]$. Definitions (3.1) and (3.14) require us to compare then, in particular, $B_{\theta_l^{(n)}}$ and $B_{\theta_k^{(n)}}$, provided $lTn^{-1}\le t<(l+1)Tn^{-1}$ and $\theta_k^{(n)}\le t<\theta_{k+1}^{(n)}$. Via standard renewal theory arguments, we conclude that in average $|k-l|$ for such $k,l\le n$ is of order $n^{1/2}$, then $|\theta_k^{(n)}-\theta_l^{(n)}|$ is of order $n^{-1/2}$, and so $|B_{\theta_k^{(n)}}-B_{\theta_l^{(n)}}|$ is roughly of order $n^{-1/4}$. The proof of the following result in the next section makes these heuristic arguments precise and an effort is made to obtain as best as possible estimates here.

Lemma 3.2. *There exists a constant $C>0$ such that, for all $n,z>0$,*

$$\begin{aligned}&|V^{B,n}(z)-V^{B,\theta,n}(z)|\\ &\qquad\le C(F_0(z)+\Delta_0(z)+z+1)n^{-1/4}(\ln n)^{3/4}.\end{aligned}\tag{3.17}$$

The values $V^{B,\theta,n}(z)$ are still defined for piecewise constant approximations $S^{B,\theta,n}(z)$ of the BS stock prices $S^B$ given by (2.4). Thus, on the next step we replace $S^{B,\theta,n}(z)$ by $S^B$ estimating the corresponding error which turns out to be of smaller order than in other lemmas, as we have to compare the Brownian motion here at times $s,t$ such that $|t-s|\le\theta_k^{(n)}-\theta_{k-1}^{(n)}$ and since the latter is of order $1/n$, the increment $|B_t-B_s|$ is roughly of order $n^{-1/2}$ which is made precise in the following result.

Lemma 3.3. *For each $\varepsilon>0$, there exists $C_\varepsilon>0$ such that, for all $z,n>0$ and $\zeta,\eta\in\mathcal{T}_{0,n}^{B,n}$,*

$$\begin{aligned}&E^B|Q_z^{B,\theta,n}(\theta_\zeta^{(n)},\theta_\eta^{(n)})-Q_z^B(\theta_\zeta^{(n)},\theta_\eta^{(n)})|\\ &\qquad\le E^B\max_{0\le k,l\le n}|Q_z^{B,\theta,n}(\theta_k^{(n)},\theta_l^{(n)})-Q_z^B(\theta_k^{(n)},\theta_l^{(n)})|\\ &\qquad\le C_\varepsilon z n^{\varepsilon-1/2}.\end{aligned}\tag{3.18}$$

In (2.10) for $V(z)$ the allowed stopping times take values in the interval $[0,T]$, so we have to restrict the stopping times $\theta_k^{(n)}$ (which are not bounded) to this interval. It is not difficult to understand that in the average $|\theta_n^{(n)}-T|$ is of order $n^{-1/2}$ and $\theta_n^{(n)}-\theta_n^{(n)}\wedge T$ is of the same order. Then the absolute value of the increment of the Brownian motion taken at times $\theta_n^{(n)}\wedge T$ and $\theta_n^{(n)}$ is roughly of order $n^{-1/4}$, and so the restriction of embedding times to the interval $[0,T]$ leads to a difference of about that order (see also Remark 3.7 below).

Lemma 3.4. *There exists a constant $C>0$ such that, for all $z,n>0$*

*and* $\zeta,\eta\in\mathcal{T}_{0,n}^{B,n}$,

$$\begin{aligned} & E^B|Q_z^B(\theta_\zeta^{(n)},\theta_\eta^{(n)})-Q_z^B(\theta_\zeta^{(n)}\wedge T,\theta_\eta^{(n)}\wedge T)| \\ & \quad\le E^B\max_{0\le k,l\le n}|Q_z^B(\theta_k^{(n)},\theta_l^{(n)})-Q_z^B(\theta_k^{(n)}\wedge T,\theta_l^{(n)}\wedge T)| \\ & \quad\le C(F_0(z)+\Delta_0(z)+z+1)n^{-1/4}. \end{aligned} \tag{3.19}$$

Until now we considered only stopping times $\theta_k^{(n)}$ for $k=0,1,\dots,n$, which may not be enough, in principle, in order to approximate all Brownian stopping times bounded by $T$, so the next result asserts that we can employ the whole sequence $\theta_0^{(n)}=0,\theta_1^{(n)},\theta_2^{(n)},\dots$. The estimates of the corresponding difference here are similar to Lemma 3.4 and they produce, essentially, the same result.

Lemma 3.5. *There exists a constant $C>0$ such that, for all $z,n>0$ and $\zeta,\eta\in\mathcal{T}^{B,n}$ (with $\mathcal{T}^{B,n}$ defined at the beginning of Section 3),*

$$\begin{aligned} & E^B|Q_z^B(\theta_\zeta^{(n)}\wedge T,\theta_\eta^{(n)}\wedge T)-Q_z^B(\theta_{\zeta\wedge n}^{(n)}\wedge T,\theta_{\eta\wedge n}^{(n)}\wedge T)| \\ & \quad\le E^B\sup_{0\le k,l<\infty}|Q_z^B(\theta_k^{(n)}\wedge T,\theta_l^{(n)}\wedge T)-Q_z^B(\theta_{k\wedge n}^{(n)}\wedge T,\theta_{l\wedge n}^{(n)}\wedge T)| \\ & \quad\le C(F_0(z)+\Delta_0(z)+z+1)n^{-1/4}. \end{aligned} \tag{3.20}$$

Set $\mathcal{T}_T^{B,n}=\{\theta_\zeta^{(n)}\wedge T:\zeta\in\mathcal{T}^{B,n}\}$ and let

$$V_{0,T}^{B,n}(z)=\inf_{\sigma\in\mathcal{T}_T^{B,n}}\sup_{\tau\in\mathcal{T}_T^{B,n}}E^BQ_z^B(\sigma,\tau). \tag{3.21}$$

Then Lemmas 3.3–3.5 yield that, for some constant $C>0$,

$$\begin{aligned} & |V^{B,\theta,n}(z)-V_{0,T}^{B,n}(z)| \\ & \quad\le\sup_{\zeta\in\mathcal{T}^{B,n}}\sup_{\eta\in\mathcal{T}^{B,n}}E^B(|Q_z^{B,\theta,n}(\theta_{\zeta\wedge n}^{(n)},\theta_{\eta\wedge n}^{(n)})-Q_z^B(\theta_{\zeta\wedge n}^{(n)},\theta_{\eta\wedge n}^{(n)})| \\ & \qquad+|Q_z^B(\theta_{\zeta\wedge n}^{(n)},\theta_{\eta\wedge n}^{(n)})-Q_z^B(\theta_{\zeta\wedge n}^{(n)}\wedge T,\theta_{\eta\wedge n}^{(n)}\wedge T)| \\ & \qquad+|Q_z^B(\theta_{\zeta\wedge n}^{(n)}\wedge T,\theta_{\eta\wedge n}^{(n)}\wedge T) \\ & \qquad\qquad-Q_z^B(\theta_\zeta^{(n)}\wedge T,\theta_\eta^{(n)}\wedge T)|) \\ & \quad\le 3C(F_0(z)+\Delta_0(z)+z+1)n^{-1/4}. \end{aligned} \tag{3.22}$$

In definition (3.21) of $V_{0,T}^{B,n}(z)$ we consider only stopping times of the special form, while in (2.11), which gives $V(z)$, all Brownian stopping times with values in $[0,T]$ are allowed and the last step in the proof of Theorem 2.1 is to estimate the corresponding error which turns out to be of the same order as in Lemma 3.3 since, again, we have to estimate increments $|B_t - B_s|$ when $|t-s| \le \theta_k^{(n)} - \theta_{k-1}^{(n)}$, though here $k$ runs over all positive integers and not only up to $n$ as in Lemma 3.3.

LEMMA 3.6. *For any $\varepsilon > 0$, there exists a constant $C_\varepsilon > 0$ such that, for all $z, n > 0$,*

$$|V(z) - V_{0,T}^{B,n}(z)| \le C_\varepsilon (F_0(z) + \Delta_0(z) + z + 1) n^{\varepsilon - 1/2}. \tag{3.23}$$

Lemmas 3.1–3.6 yield the required estimate of $|V(z) - V^{(n)}(z)|$ from Theorem 2.1 and the corresponding estimate of $|V(z) - \hat{V}^{(n)}(z)|$ goes through exactly in the same way.

REMARK 3.7. The estimate of Theorem 2.1 (and so the estimates of Theorems 2.2 and 2.3) seems to be, essentially, optimal under the general conditions (2.1) and (2.2) at least, using the method which relies on the Skorokhod embedding as above. It is known and can be seen from the proof of Lemma 3.2 that the embedding procedure cannot provide, in general, a better than $n^{-1/4}$ estimate there. One may restrict the class of payoffs assuming, for instance, that for piecewise constant functions $v$ of time $t \in [0, T]$, the functionals $F_t(v)$ and $\Delta_t(v)$ depend only on the values of $v$ but not on the time intervals between jumps of $v$. This condition is satisfied, for instance, in the case of Russian type options. Then we can skip Lemma 3.2 and after a slight modification of Lemma 3.1, we can proceed directly to Lemma 3.3. In view of Lemmas 3.3–3.6, this would lead to a slightly better estimate $Cn^{-1/4}$ than the estimate (2.13) of Theorem 2.1. Still, it does not seem possible to obtain under reasonably general conditions (which are satisfied, say, for Russian options) better that $n^{-1/4}$ estimates in Lemmas 3.4 and 3.5. Indeed, in order to obtain specific estimates, we have to get rid of the general functionals $F$ and $\Delta$ using the assumptions (2.1) and (2.2), which inevitably leads to an estimate of

$$E^B \sup_{\theta_n^{(n)} \wedge T \le t \le \theta_n^{(n)} \vee T} |B_t - B_{\theta_k^{(n)} \wedge T}|$$

(in fact, of a bit larger expression), which by the Burkholder–Davis–Gandy inequality (see [11], Theorem 3.1 in Section 3.3 and [13], Theorem 3.28 in Section 3.3) is of order

$$E^B (\theta_n^{(n)} \vee T - \theta_n^{(n)} \wedge T)^{1/2}$$

and the latter expression is of order $n^{-1/4}$. The main obstruction to a better estimate of the first expression above is that (2.2) requires us to write the supremum and the absolute value inside and not outside of the expectation. This obstruction disappears in the other two papers [24] and [30] employing the Skorokhod embedding which also have to face estimates of the error originated from the fact that, after embedding, we have to consider the Brownian motion until the stopping time $\theta_n^{(n)}$ which differs from the expiry time $T$ by about $n^{-1/2}$. In [30] only European options with payoffs depending on the current stock price are considered, which enables the author to apply the simple random walk machinery leading to better estimates. In [24] and in its generalization considered in Section 6 below the payoffs also depend only on the current stock price, which together with the smoothness assumption enables us to use the Itô formula leading to the Dynkin formula where the stochastic integral part disappears and the remaining Riemann integral taken between $T$ and $\theta_n^{(n)}$ has the same order $n^{-1/2}$ as $|T-\theta_n^{(n)}|$. There exist other methods of strong invariance principle type uniform approximations of the Brownian motion by means of properly normalized sums of i.i.d. random variables (see, e.g., [3, 16, 31]) which may give a better rate of approximation, but the problem arising there is to find a proper substitution to Lemma 3.1 which would enable us not only to consider corresponding processes on one probability space, but also to deal with stopping times with respect to the same filtration in the inf sup formulas expressing values of corresponding Dynkin's games. In the case of European options (or contingent claims) with path dependent payoffs satisfying (2.1), we do not have to worry about stopping times and need only to produce a best possible uniform approximation of BS stock prices by appropriate CRR stock prices. Employing the quantile transformation method from [16] and [31], this can always be done with an error (roughly) of order $n^{-1/2}$. On the other hand, the method of [3] can be used, in principle, in order to approximate markets where stock prices evolve not necessarily as a geometric Brownian motion.

Remark 3.8. It follows from [14] that, with probability one,

$$\limsup_{n\to\infty}(|B_{\theta_n^{(n)}}-B_T|n^{1/4}(\ln n)^{-1/2}(\ln\ln n)^{-1/4})=2^{1/4}.$$

It would be interesting to understand whether the estimate (2.13) of Theorem 2.1 can be improved to $Cn^{-1/4}\sqrt{\ln n}(\ln\ln n)^{1/4}$ or the present estimate is sharp. In view of Lemma 4.1 below, the estimate (2.13) is connected with the bound for

$$E^B\max_{0\le k\le n}|B_{\theta_k^{(n)}}-B_{kT/n}|,$$

though (2.13) requires also an estimate of $H_2^{(n)}$ in the next section which does not seem to be reducible to this one. This question can be formulated in the following classical form. Let $\hat{\Theta}_0 = 0$ and, successively, $\hat{\Theta}_{n+1} = \inf\{t > \hat{\Theta}_n : |B_t - B_{\hat{\Theta}_n}| = 1\}$ with $B_0 = 0$. The result of [14] gives that, with probability one,

$$\limsup_{n\to\infty}(|B_{\hat{\Theta}_n} - B_n| n^{-1/4}(\ln n)^{-1/2}(\ln\ln n)^{-1/4}) = 2^{1/4}.$$

For our problem, we need to know the asymptotical behavior as $n \to \infty$ of

$$E^B \max_{0\le k\le n} |B_{\hat{\Theta}_k} - B_k|.$$

Our estimates give the bound $Cn^{1/4}(\ln n)^{3/4}$ for this expectation. Is there a better bound or this bound the best possible?

## 4. Proving the estimates.

Set

$$B_t^{(n)} = -\frac{\kappa t}{2}\sqrt{\frac{T}{n}} + B_t \quad \text{and} \quad \Theta^{(n)} = \inf\{t > 0 : |B_t^{(n)}| = 1\}.$$

By the scaling property of the Brownian motion,

$$\sqrt{\frac{T}{n}} B_t^{(n)} \overset{\mathrm{d}}{\sim} B^*_{(T/nt)} \quad \text{and} \quad \theta_1^{(n)} \overset{\mathrm{d}}{\sim} \frac{T}{n}\Theta^{(n)}, \tag{4.1}$$

where $\xi \overset{\mathrm{d}}{\sim} \tilde{\xi}$ means that $\xi$ and $\tilde{\xi}$ have the same distribution. Observe that, in view of independency of increments $B_l - B_{l-1}, l = 1, 2, \dots$, for any $n \ge 1$,

$$P^B\{\Theta^{(n)} \ge k\} \le P^B\{|B_l - B_{l-1}| \le 2 + \kappa\sqrt{T}\ \forall\, l = 1, \dots, k\} = e^{-b_T k}, \tag{4.2}$$

where $b_T = -\ln P^B\{|B_1| \le 2 + \kappa\sqrt{T}\} > 0$. It follows that, for any nonnegative $a < b_T$,

$$E^B e^{a\Theta^{(n)}} \le \sum_{k=0}^{\infty} e^{a(k+1)} P^B\{\Theta^{(n)} \ge k\} \le e^a(1 - e^{a-b_T})^{-1} < \infty. \tag{4.3}$$

The estimates of Section 3 are not trivial only for large $n$, so we will assume that $n$ is sufficiently big, in particular, that various exponential moments of the form $E^B \exp(\frac{aT}{n}\Theta^{(n)})$ are finite, that is, that $n > aTb_T^{-1}$.

Since $\exp \kappa B_t^*$, $t \ge 0$, is a martingale with respect to the probability $P^B$, and so $E^B \exp(\kappa B^*_{\theta_1^{(n)}}) = 1$ (assuming that $n > \frac{1}{2}Tb_T^{-1}$), we derive by an easy computation that $B^*_{\theta_1^{(n)}} = \sqrt{\frac{T}{n}}$ or $= -\sqrt{\frac{T}{n}}$ and $B^{(n)}_{\Theta^{(n)}} = 1$ or $= -1$ with probability $(1 + \exp(\kappa\sqrt{\frac{T}{n}}))^{-1}$ or $(1 + \exp(-\kappa\sqrt{\frac{T}{n}}))^{-1}$, respectively. Set $\alpha_n = E^B\Theta^{(n)}$ so that $E^B\theta_1^{(n)} = \alpha_n\frac{T}{n}$. Since the Brownian motion is a martingale, and so $E^B B_{\Theta^{(n)}} = 0$, we have that

$$-\frac{\kappa}{2}\alpha_n\sqrt{\frac{T}{n}} = E^B\left(-\frac{\kappa}{2}\Theta^{(n)}\sqrt{\frac{T}{n}} + B_{\Theta^{(n)}}\right)$$
$$= (1+e^{\kappa\sqrt{T/n}})^{-1} - (1+e^{-\kappa\sqrt{T/n}})^{-1}.$$

This together with an easy estimate shows that

$$|\alpha_n - 1| \le \min\left(2\kappa^{-1}\sqrt{\frac{n}{T}}, \frac{\kappa^2 T}{2n}\left|1-\frac{\kappa^2 T}{n}\right|^{-1}\right) \le K_1\frac{T}{n}, \tag{4.4}$$

where

$$K_1 = \min(2\kappa^{-1}(2\kappa^2 + T^{-1})^{3/2}, \kappa^2).$$

By (4.2),

$$\begin{aligned} E^B|\Theta^{(n)}|^m &\le M_m = \sum_{k=1}^{\infty} k^m e^{-b_T(k-1)} \\ &\le e^{2b_T}\int_0^\infty x^m e^{-b_T x}\,dx = e^{2b_T} m! b_T^{-(m+1)}. \end{aligned} \tag{4.5}$$

Assuming, without loss of generality, that $n \ge K_1 T$, we obtain from (4.4) that

$$P^B\{|\Theta^{(n)} - \alpha_n|^m \ge k\} \le P^B\{\Theta^{(n)} \ge k-2\},$$

and so

$$E^B|\Theta^{(n)} - \alpha_n|^m \le \sum_{k=1}^{\infty} k^m e^{-b_T(k-3)} = e^{2b_T} M_m. \tag{4.6}$$

Observe that $\theta_k^{(n)} - \alpha_n k\frac{T}{n}$, $k = 0,1,2,\dots$, is a martingale with respect to the filtration $\mathcal{F}^B_{\theta_k^{(n)}}$, $k \in \mathbb{N}$. Thus, using that $(a+b)^m \le 2^{m-1}(a^m + b^m)$ for $a, b \ge 0$, $m \ge 1$, we obtain by (4.1), (4.4), (4.6) and the Burkholder–Davis–Gandy inequality (see [11], Theorem 3.1 in Section 3.3 and [13], Theorem 3.28 in Section 3.3) that, for any $m > 1/2$,

$$\begin{aligned} E^B &\sup_{\zeta\in\mathcal{T}^{B,n}_{0,n}} \left|\theta^{(n)}_\zeta - \zeta\frac{T}{n}\right|^{2m} \\ &= E^B \max_{0\le k\le n}\left|\theta_k^{(n)} - \frac{kT}{n}\right|^{2m} \\ &\le 2^{2m-1}T^{2m}|\alpha_n - 1|^{2m} + 2^{2m-1}E^B\max_{0\le k\le n}\left|\theta_k^{(n)} - \alpha_n\frac{kT}{n}\right|^{2m} \end{aligned}$$

$$
\begin{aligned}
&\le 2^{2m-1}K_1^{2m}T^{4m}n^{-2m} + 2^{2m-1}\Lambda_m\left(nE^B\left|\theta_1^{(n)} - \alpha_n\frac{T}{n}\right|^2\right)^m \\
(4.7)\qquad &\le 2^{2m-1}T^{2m}n^{-m}(T^{2m}K_1^{2m}n^{-m} + \Lambda_m e^{2mb_T}M_2^m) \\
&\le K_2^{(m)}T^{2m}n^{-m},
\end{aligned}
$$

where $\Lambda_m = 4^{m^2}m^{m(2m+1)}(2m-1)^{m(1-2m)}$ and $K_2^{(m)} = 2^{2m-1}(K_1^{2m}T^{2m} + \Lambda_m e^{2mb_T}M_2^m)$ assuming that $n \ge 1$. We will need (4.7) mostly with $m = 1$ which requires only the Doob–Kolmogorov inequality (see, e.g., [11]) and in order to simplify notation, we set $K_2 = K_2^{(1)}$.

Using the exponential martingale $\exp(aB_t - \frac{1}{2}a^2t)$ and applying the Doob–Kolmogorov and Cauchy–Schwarz inequalities, we obtain

$$
\begin{aligned}
E^B \sup_{0\le t\le\tau} \exp(aB_t) &\le E^B e^{(1/2)a^2\tau} \sup_{0\le t\le\tau} \exp(aB_t - \tfrac{1}{2}a^2t) \\
&\le (E^B e^{a^2\tau})^{1/2}\left(E^B \sup_{0\le t\le\tau}\exp(2aB_t - a^2t)\right)^{1/2} \\
(4.8)\qquad &\le 2(E^B e^{a^2\tau})^{1/2}(E^B\exp(2aB_\tau - a^2\tau))^{1/2} \\
&\le 2(E^B e^{a^2\tau})^{1/2}(E^B\exp(4aB_\tau - 8a^2\tau))^{1/4}(E^B e^{6a^2\tau})^{1/4} \\
&= 2(E^B e^{a^2\tau})^{1/2}(E^B e^{6a^2\tau})^{1/4}
\end{aligned}
$$

for any finite Brownian stopping time $\tau$ and a number $a$. If $\sigma \le \tau$ is another Brownian stopping time, then by the Burkholder–Davis–Gandy inequality (see [11] and [13]) applied to the (continuous) martingale (stochastic integral) $\int_0^t \mathbb{I}_{\sigma<s\le\tau}\,dB_s$ we obtain that, for any $m > 0$,

$$
(4.9)\qquad E^B \sup_{\sigma\le t\le\tau} |B_t - B_\sigma|^{2m} \le \Lambda_m E^B|\tau - \sigma|^m,
$$

where $\Lambda_m$ is the same as in (4.7) and, again, we will use (4.9) only for $m > 1/2$. Recall that our relevant formulas do not depend on a particular choice of a continuous in time version of the Brownian motion $B_t$ and each such version is, in fact, Hölder continuous with probability one.

In the proof of Lemma 3.2 we will need also certain renewal theory estimates which seem to be standard, but, since we could not find a direct reference, their proof for readers' convenience is given here.

LEMMA 4.1. *Let* $k_t^{(n)} = \max\{j \le n : \theta_j^{(n)} \le t\}$ *for all* $t \ge 0$ *and* $\ell_t^{(n)} = [nt/T]$ *if* $t \in [0,T]$ *and* $\ell_t^{(n)} = n$ *if* $t > T$. *Then*

$$
(4.10)\qquad E^B \sup_{0\le t\le T}|k_t^{(n)} - \ell_t^{(n)}|^2 \le 2E^B \sup_{0\le t\le T}|k_t^{(n)} - nt/T|^2 + 2 \le 2(K_2+2)n
$$

*and*

$$E^B \sup_{0\le t\le \theta_n^{(n)}} |B_{\theta^{(n)}_{k_t^{(n)}}} - B_{\theta^{(n)}_{\ell_t^{(n)}}}|^2 \le K_3 n^{-1/2}(\ln n)^{3/2}, \tag{4.11}$$

*where* $K_3>0$ *can be estimated from the proof below.*

PROOF. Let $\Theta_1^{(n)},\Theta_2^{(n)},\dots$ be i.i.d. random variables with the same distribution as $\Theta^{(n)}$. Set $m_u^{(n)}=\max\{j\le n:\sum_{i=1}^j\Theta_i^{(n)}\le u\}$, then, by (4.1), the process $m_{nt/T}^{(n)}$, $t\in[0,T]$, has the same distribution as the process $k_t^{(n)}$, $t\in[0,T]$. Hence,

$$E^B\sup_{0\le t\le T}|k_t^{(n)}-nt/T|^2 = E^B\sup_{0\le u\le n}|m_u^{(n)}-u|^2.$$

Set $\Psi_t=\sum_{j=1}^{[t]}\Theta_j^{(n)}$ for $t\ge 1$ and $\Psi_t=\Psi_0=0$ for $t\in[0,1)$. It is clear that if $l<n$, then $m_u^{(n)}-u=l-u$ if and only if $l-\Psi_l\ge l-u>l-\Psi_{l+1}$, and so in this case

$$|m_u^{(n)}-u|\le\max(|\Psi_l-l|,|\Psi_{l+1}-(l+1)|+1).$$

If $m_u^{(n)}=n$ and $u\le n$, then $\Psi_n\le u\le n$, and so $|m_u^{(n)}-u|\le|\Psi_n-n|$. Hence,

$$\max_{0\le u\le n}|m_u^{(n)}-u|\le\max_{0\le l\le n}|\Psi_l-l|+1. \tag{4.12}$$

Observe that, by (4.4), for any $l\le n$,

$$|\Psi_l-l|\le|\Psi_l-l\alpha_n|+K_1T, \tag{4.13}$$

and so by the Doob–Kolmogorov inequality,

$$\begin{aligned}E^B\max_{0\le l\le n}|\Psi_l-l|^2 &\le 2E^B\max_{0\le l\le n}|\Psi_l-l\alpha_n|^2+2K_1^2T^2\\ &\le 8E^B|\Psi_n-n\alpha_n|^2+2K_1^2T^2=8nE^B(\Theta^{(n)}-\alpha_n)^2+2K_1^2T^2\end{aligned}$$

and (4.10) follows from (4.4) and (4.6).

Next, we prove (4.11) estimating, first,

$$\begin{aligned}&\sup_{0\le t\le\theta_n^{(n)}}|B_{\theta^{(n)}_{k_t^{(n)}}}-B_{\theta^{(n)}_{\ell_t^{(n)}}}|^2\\ &\quad\le 4\mathbb{I}_{\sup_{0\le t\le\theta_n^{(n)}}|k_t^{(n)}-\ell_t^{(n)}|>D\sqrt{n\ln n}}\sup_{0\le t\le\theta_n^{(n)}}|B_t|^2\\ &\qquad+4\mathbb{I}_{\max_{k,l\le n,|k-l|\le D\sqrt{n\ln n}}|\theta_k^{(n)}-\theta_l^{(n)}|>D^2\sqrt{n^{-1}\ln n}}\sup_{0\le t\le\theta_n^{(n)}}|B_t|^2\\ &\qquad+4\mathbb{I}_{\max_{k\le n}\sup_{0\le t\le D^2\sqrt{n^{-1}\ln n}}|B_{\theta_k^{(n)}+t}-B_{\theta_k^{(n)}}|>D^2n^{-1/4}(\ln n)^{3/4}}\end{aligned} \tag{4.14}$$

$$\times \sup_{0\le t\le\theta_n^{(n)}} |B_t|^2$$
$$+ D^4 n^{-1/2}(\ln n)^{3/2},$$

where $D>0$ will be chosen below. Observe that

$$\sup_{0\le t\le\theta_n^{(n)}} |k_t^{(n)}-\ell_t^{(n)}| \le 2 \sup_{0\le t\le T}|k_t^{(n)}-\ell_t^{(n)}| \tag{4.15}$$

since, by the definition of $k_t^{(n)}$ and $\ell_t^{(n)}$,

$$\sup_{T\le t\le T\vee\theta_n^{(n)}} |k_t^{(n)}-\ell_t^{(n)}| \le n-k_T^{(n)} = \ell_T^{(n)}-k_T^{(n)}.$$

Since $|k_t^{(n)}-\ell_t^{(n)}|\le|k_t^{(n)}-nt/T|+1$ and the processes $m_{nt/T}^{(n)}, t\in[0,T]$, and $k_t^{(n)}$, $t\in[0,T]$, have the same distribution, we derive from (4.12), (4.13) and (4.15) that

$$\begin{aligned} &P^B\Big\{\sup_{0\le t\le\theta_n^{(n)}}|k_t^{(n)}-\ell_t^{(n)}|>D\sqrt{n\ln n}\Big\}\\ &\quad\le P^B\Big\{\max_{0\le l\le n}|\Psi_l-l\alpha_n|>\tfrac12 D\sqrt{n\ln n}-2\Big\}\\ &\quad\le \sum_{l=0}^{n}(P^B\{\Psi_l-l\alpha_n>\tfrac12 D\sqrt{n\ln n}-2\}\\ &\qquad\qquad + P^B\{l\alpha_n-\Psi_l>\tfrac12 D\sqrt{n\ln n}-2\}). \end{aligned} \tag{4.16}$$

By (4.5), (4.6), Chebyshev's inequality and the definition of $\Psi_l$,

$$\begin{aligned} &P^B\Big\{\Psi_l-l\alpha_n>\frac12 D\sqrt{n\ln n}-2\Big\}\\ &\quad\le P^B\Big\{\exp\Big(2\sqrt{n^{-1}\ln n}\sum_{i=1}^{l}(\Theta_i^{(n)}-\alpha_n)\Big)\ge n^D e^{-4}\Big\}\\ &\quad\le e^4 n^{-D}(E^B\exp(2\sqrt{n^{-1}\ln n}(\Theta_i^{(n)}-\alpha_n)))^l\\ &\quad\le e^4 n^{-D}\Big(1+e^{2b_T}\sum_{m=2}^{\infty}\Big(\frac{4\ln n}{n}\Big)^{m/2}\frac{M_m}{m!}\Big)^n\\ &\quad\le e^4 n^{(8e^{4b_T}b_T^{-3}-D)}\le e^4 n^{-2}, \end{aligned} \tag{4.17}$$

where we use the inequality $(1+a/q)^q<e^a$ for $a,q>0$, choose $D\ge 2+8e^{4b_T}b_T^{-3}$ and assume that $n\ge 2(16/b_T^2)^2$, so that $n^{-1}\ln n\le b_T^2/16$. Similarly, under the same conditions,

$$P^B\{l\alpha_n-\Psi_l>\tfrac12 D\sqrt{n\ln n}-2\}\le e^4 n^{-2}. \tag{4.18}$$

Next, by (4.1), (4.3) and the Chebyshev inequality,

$$\begin{aligned}
&P^B\Big\{\max_{k,l\le n,|k-l|\le D\sqrt{n\ln n}}|\theta_k^{(n)}-\theta_l^{(n)}|>D^2\sqrt{n^{-1}\ln n}\Big\}\\
&\quad\le nP^B\{\theta^{(n)}_{[D\sqrt{n\ln n}]}>D^2\sqrt{n^{-1}\ln n}\}\\
&\quad= nP^B\{anT^{-1}\theta^{(n)}_{[D\sqrt{n\ln n}]}>aD^2T^{-1}\sqrt{n\ln n}\}\\
&\quad\le n\exp(-aD^2T^{-1}\sqrt{n\ln n})E^B\exp(anT^{-1}\theta^{(n)}_{[D\sqrt{n\ln n}]})\\
&\quad\le n(\exp(-aDT^{-1})E^Be^{a\Theta^{(n)}})^{D\sqrt{n\ln n}}\le n^{-1},
\end{aligned}\tag{4.19}$$

(where $[b]$ is the integral part of $b$) if we choose a positive $a<b_T$,

$$D\ge Ta^{-1}(1-\ln(1-e^{a-b_T}))+T+1 \qquad \text{so that } (e^{a(DT^{-1}-1)}(1-e^{a-b_T}))^D\ge e,$$

and assume that $n\ge e^3$, so that $n^{-1}\ln n\le 1/4$. Now, by the strong Markov property, the reflection principle and the scaling property of the Brownian motion (see, e.g., [13], Chapter 2),

$$\begin{aligned}
&P^B\Big\{\max_{k\le n}\sup_{\theta_k^{(n)}\le t\le\theta_k^{(n)}+D^2\sqrt{n^{-1}\ln n}}|B_t-B_{\theta_k^{(n)}}|>D^2n^{-1/4}(\ln n)^{3/4}\Big\}\\
&\quad\le nP^B\Big\{\sup_{0\le t\le D^2\sqrt{n^{-1}\ln n}}|B_t|>D^2n^{-1/4}(\ln n)^{3/4}\Big\}\\
&\quad\le 4nP^B\{B_{D^2\sqrt{n^{-1}\ln n}}>D^2n^{-1/4}(\ln n)^{3/4}\}\\
&\quad= 4nP^B\{B_1>D\sqrt{\ln n}\}\\
&\quad= 4n\int_{D(\ln n)^{1/2}}^{\infty}(2\pi)^{-1/2}e^{-x^2/2}\,dx\\
&\quad\le\frac{8n}{D\sqrt{2\pi\ln n}}e^{-(1/2)D^2\ln n}\le\frac{4}{n\sqrt{2\pi\ln n}},
\end{aligned}\tag{4.20}$$

provided we choose $D\ge 2$. Finally, by (4.1), (4.4) and (4.9),

$$E^B\sup_{0\le t\le\theta_n^{(n)}}|B_t|^2\le 4E^B\theta_n^{(n)}=4T\alpha_n\le 4T(1+K_1T),\tag{4.21}$$

and we obtain (4.11) from (4.14) and (4.16)–(4.21) together with the Cauchy–Schwarz inequality. □

Now we are ready to pass directly to the proof of Lemma 3.2.

Proof of Lemma 3.2. By (3.3), (3.16) and the equality $\theta_\zeta^{(n)} \wedge \theta_\eta^{(n)} = \theta_{\zeta\wedge\eta}^{(n)}$,

$$|V^{B,n}(z) - V^{B,\theta,n}(z)| \le \sup_{\zeta\in\mathcal{T}_{0,n}^{B,n}} \sup_{\eta\in\mathcal{T}_{0,n}^{B,n}} (J_1(\zeta,\eta) + J_2(\zeta,\eta)), \tag{4.22}$$

where

$$J_1(\zeta,\eta) = E^B\left(|e^{-rT/n\zeta\wedge\eta} - e^{-r\theta_{\zeta\wedge\eta}^{(n)}}|R_z^{B,n}\left(\frac{\zeta T}{n},\frac{\eta T}{n}\right)\right)$$

and

$$J_2(\zeta,\eta) = E^B\left|R_z^{B,n}\left(\frac{\zeta T}{n},\frac{\eta T}{n}\right) - R_z^{B,\theta,n}(\theta_\zeta^{(n)},\theta_\eta^{(n)})\right|.$$

Since $|e^{-ra} - e^{-rb}| \le r|a-b|$, we obtain by (4.7) and the Cauchy–Schwarz inequality

$$\begin{aligned} J_1(\zeta,\eta) &\le rE^B\left|\theta_{\zeta\wedge\eta}^{(n)} - \frac{T}{n}\zeta\wedge\eta\right|R_z^{B,n}\left(\frac{\zeta T}{n},\frac{\eta T}{n}\right) \\ &\le \frac{\sqrt{K_2}rT}{\sqrt{n}}(J_{11}(\eta,\zeta))^{1/2}, \end{aligned} \tag{4.23}$$

where by (2.3), (3.1), (3.2) and (4.8),

$$\begin{aligned} J_{11}(\zeta,\eta) &= E^B\left(R_z^{B,n}\left(\frac{\zeta T}{n},\frac{\eta T}{n}\right)\right)^2 \\ &\le 2E^B((F_{(T/n)\zeta\wedge\eta}(S^{B,n}(z)))^2 + (\Delta_{(T/n)\zeta\wedge\eta}(S^{B,n}(z)))^2) \\ &\le 6\left(F_0^2(z) + \Delta_0^2(z) + 2L^2(T+2)^2E^B\left(1 + \sup_{0\le t\le T}(S_t^{B,n}(z))^2\right)\right) \\ &\le 6\left(F_0^2(z) + \Delta_0^2(z) + 2L^2(T+2)^2\left(1 + z^2e^{2rT}E^B \sup_{0\le t\le\theta_n^{(n)}} e^{2\kappa B_t}\right)\right) \\ &\le 6(F_0^2(z) + \Delta_0^2(z) \\ &\qquad + 4L^2(T+2)^2(1 + z^2e^{2rT}(E^Be^{4\kappa^2\theta_n^{(n)}})^{1/2}(E^Be^{24\kappa^2\theta_n^{(n)}})^{1/4})). \end{aligned} \tag{4.24}$$

Since

$$\theta_n^{(n)} = \sum_{k=1}^n(\theta_k^{(n)} - \theta_{k-1}^{(n)}) \quad \text{and} \quad B_{\theta_n^{(n)}} = \sum_{k=1}^n(B_{\theta_k^{(n)}} - B_{\theta_{k-1}^{(n)}})$$

are sums of i.i.d. random variables, we obtain by (4.1)–(4.5) and the Taylor formula that, for any $a > 0$,

$$E^Be^{a\theta_n^{(n)}} = (E^Be^{a(T/n)\Theta^{(n)}})^n$$

$$
\begin{aligned}
&\leq \left(1 + a(K_1+1)\frac{T}{n} + \frac{a^2T^2}{n^2}\sum_{m=0}^{\infty}\frac{a^mT^mM_{m+2}}{n^m(m+2)!}\right)^n \\
&\leq \left(1 + a(K_1+1)\frac{T}{n} + a^2T^2b^{-3}n^{-2}\exp\left(2b_T + \frac{aT}{nb_T}\right)\right)^n \\
&\leq C_a = e^{a(K_1+1)T},
\end{aligned} \tag{4.25}
$$

provided $n \geq Tb_T^{-3}\exp(2b_T + aTb_T^{-1})$ and we use that $(1+a/q)^q \leq e^a$ if $a, q > 0$. This together with (4.8) gives also

$$
E^B e^{aB_{\theta_n^{(n)}}} \leq 2(E^B e^{a^2\theta_n^{(n)}})^{1/2}(E^B e^{6a^2\theta_n^{(n)}})^{1/4} \leq 2(C_{a^2})^{1/2}(C_{6a^2})^{1/4}, \tag{4.26}
$$

assuming that $n \geq Tb_T^{-3}\exp(2b_T + 2a^2Tb_T^{-1})$.

Next, we estimate $J_2(\zeta,\eta)$. By (3.2) and (3.15),

$$
\begin{aligned}
J_2(\zeta,\eta) \leq E^B(&|F_{\zeta T/n}(S^{B,n}(z)) - F_{\theta_\zeta^{(n)}}(S^{B,\theta,n}(z))| \\
&+ |F_{\eta T/n}(S^{B,n}(z)) - F_{\theta_\eta^{(n)}}(S^{B,\theta,n}(z))| \\
&+ |\Delta_{\zeta T/n}(S^{B,n}(z)) - \Delta_{\theta_\zeta^{(n)}}(S^{B,\theta,n}(z))|).
\end{aligned} \tag{4.27}
$$

For any $\zeta \in \mathcal{T}_{0,n}^{B,n}$, we obtain from (2.1) and (2.2) that

$$
\begin{aligned}
&|F_{\zeta T/n}(S^{B,n}(z)) - F_{\theta_\zeta^{(n)}}(S^{B,\theta,n}(z))| \\
&\quad \leq \max_{0\leq k\leq n}(|F_{\theta_k^{(n)}}(S^{B,n}(z)) - F_{\theta_k^{(n)}}(S^{B,\theta,n}(z))| \\
&\qquad + |F_{kT/n}(S^{B,n}(z)) - F_{\theta_k^{(n)}}(S^{B,n}(z))|) \\
&\quad \leq L(\theta_n^{(n)}+1)J_{21} + L\max_{0\leq k\leq n}\left|\theta_k^{(n)} - \frac{kT}{n}\right|(1+J_{22}) + LJ_{23},
\end{aligned} \tag{4.28}
$$

where

$$
\begin{aligned}
J_{21} &= \sup_{0\leq t\leq\theta_n^{(n)}}|S_t^{B,n}(z) - S_t^{B,\theta,n}(z)|, \\
J_{22} &= \sup_{0\leq t\leq T}S_t^{B,n}(z) \leq ze^{rT}\sup_{0\leq t\leq\theta_n^{(n)}}e^{\kappa B_t}
\end{aligned} \tag{4.29}
$$

and

$$
J_{23} = \max_{0\leq k\leq n}\sup_{(kT/n)\wedge\theta_k^{(n)}\leq u\leq t\leq(kT/n)\vee\theta_k^{(n)}}|S_t^{B,n}(z) - S_u^{B,n}(z)|.
$$

Set

$$
H_1^{(n)}(t) = \left|r\left(\theta_{k_t^{(n)}}^{(n)} - \ell_t^{(n)}\frac{T}{n}\right) + \kappa(B^*_{\theta_{k_t^{(n)}}^{(n)}} - B^*_{\theta_{\ell_t^{(n)}}^{(n)}})\right| \quad \text{and} \quad H_1^{(n)} = \sup_{0\leq t\leq\theta_n^{(n)}}H_1^{(n)}(t),
$$

with $k_t^{(n)}$ and $\ell_t^{(n)}$ defined in Lemma 4.1. Then by (3.1), (3.14) and (4.29),

$$\begin{aligned} J_{21} &\le J_{22}\mathbb{I}_{H_1^{(n)}\le 1}\sup_{0\le t\le\theta_n^{(n)}}|e^{H_1^{(n)}(t)}-1| \\ &\quad + \mathbb{I}_{H_1^{(n)}>1}\sup_{0\le t\le\theta_n^{(n)}}(S_t^{B,n}(z)+S_t^{B,\theta,n}(z)) \\ &\le 2ze^{rT}\sup_{0\le t\le\theta_n^{(n)}}e^{\kappa B_t}H_1^{(n)}+2z\mathbb{I}_{H_1^{(n)}>1}(e^{rT}+e^{r\theta_n^{(n)}})\Big(\sup_{0\le t\le\theta_n^{(n)}}e^{\kappa B_t}\Big). \end{aligned} \tag{4.30}$$

Next, by (4.15),

$$\begin{aligned} H_1^{(n)} &\le \Big|r-\frac{\kappa^2}{2}\Big|\max_{0\le k\le n}\Big|\theta_k^{(n)}-k\frac{T}{n}\Big| \\ &\quad +2\frac{T}{n}\Big|r-\frac{\kappa^2}{2}\Big|\sup_{0\le t\le T}|k_t^{(n)}-\ell_t^{(n)}|+\kappa\sup_{0\le t\le\theta_n^{(n)}}|B_{\theta^{(n)}_{k_t^{(n)}}}-B_{\theta^{(n)}_{\ell_t^{(n)}}}|. \end{aligned} \tag{4.31}$$

Hence, by (4.7), (4.10), (4.11), (4.31), the Cauchy–Schwarz and the Chebyshev inequalities, it follows that there exists a constant $\tilde C>0$ such that

$$P^B\{H_1^{(n)}>1\}\le\tilde Cn^{-1/2}(\ln n)^{3/2}. \tag{4.32}$$

This together with (4.7), (4.8), (4.10), (4.11), (4.25), (4.30), (4.31) and the Cauchy–Schwarz inequalities yields that there exists a constant $C^{(1)}>0$ such that

$$E^BJ_{21}\le C^{(1)}zn^{-1/4}(\ln n)^{3/4} \tag{4.33}$$

and both $\tilde C$ and $C^{(1)}$ can be estimated from the above formulas.

In order to estimate $J_{23}$, set

$$H_2^{(n)}(s,t)=\frac{rT}{n}(\ell_t^{(n)}-\ell_s^{(n)})+\kappa(B^*_{\theta^{(n)}_{\ell_t^{(n)}}}-B^*_{\theta^{(n)}_{\ell_s^{(n)}}})$$

and

$$H_2^{(n)}=\max_{0\le k\le n}\sup_{(kT/n)\wedge\theta_k^{(n)}\le s\le t\le(kT/n)\vee\theta_k^{(n)}\wedge T}H_2^{(n)}(s,t).$$

Then by (3.1), similarly to (4.30),

$$J_{23}\le 2J_{22}(\mathbb{I}_{H_2^{(n)}>1}+H_2^{(n)}). \tag{4.34}$$

If $\frac{kT}{n}\wedge\theta_k^{(n)}\le s\le t\le\frac{kT}{n}\vee\theta_k^{(n)}$, then by (4.4),

$$\begin{aligned} \ell_t^{(n)}-\ell_s^{(n)} &\le\frac{n}{T}(t-s)+1\le\frac{n}{T}\max_{0\le k\le n}\Big|\theta_k^{(n)}-k\frac{T}{n}\Big|+1 \\ &\le\frac{n}{T}\max_{0\le k\le n}\Big|\theta_k^{(n)}-k\frac{T}{n}\alpha_n\Big|+K_1T+1 \end{aligned} \tag{4.35}$$

and

$$
\begin{aligned}
\theta^{(n)}_{\ell^{(n)}_t} - \theta^{(n)}_{\ell^{(n)}_s} &\le \left|\theta^{(n)}_{\ell^{(n)}_t} - \frac{\ell^{(n)}_t T}{n}\right| + \left|\theta^{(n)}_{\ell^{(n)}_s} - \frac{\ell^{(n)}_s T}{n}\right| + \frac{T}{n}|\ell^{(n)}_t - \ell^{(n)}_s| \\
&\le 3 \max_{0\le k\le n}\left|\theta^{(n)}_k - \frac{kT}{n}\right| + \frac{T}{n} \qquad (4.36)\\
&\le 3 \max_{0\le k\le n}\left|\theta^{(n)}_k - \frac{kT}{n}\alpha_n\right| + 3K_4 n^{-1},
\end{aligned}
$$

where $K_4 = T(K_1T+1)$. Hence, similarly to (4.14) and (4.31),

$$
\begin{aligned}
H^{(n)}_2 &\le \left|r - \frac{3\kappa^2}{2}\right| \max_{0\le k\le n}\left|\theta^{(n)}_k - k\frac{T}{n}\right| + \left|1 - \frac{\kappa^2}{2}\right|\frac{T}{n} \\
&\quad + \max_{0\le k\le n} \sup_{t\le\theta^{(n)}_n - \theta^{(n)}_k, 0\le t\le 3\max_{0\le l\le n}|\theta^{(n)}_l - (lT/n)\alpha_n| + 3K_4n^{-1}} |B_{\theta^{(n)}_k + t} \qquad (4.37)\\
&\qquad - B_{\theta^{(n)}_k}|
\end{aligned}
$$

and

$$
\begin{aligned}
&\sup_{k\le n, t\le\theta^{(n)}_n - \theta^{(n)}_k, 0\le t\le 3\max_{0\le l\le n}|\theta^{(n)}_l - (lT/n)\alpha_n| + 3K_4n^{-1}} |B_{\theta^{(n)}_k + t} - B_{\theta^{(n)}_k}| \\
&\quad \le 2\mathbb{I}_{\max_{0\le l\le n}|\theta^{(n)}_l - (lT/n)\alpha_n| > (1/3)D\sqrt{n^{-1}\ln n} - K_4n^{-1}} \\
&\qquad \times \sup_{0\le t\le \theta^{(n)}_n} |B_t| + Dn^{-1/4}(\ln n)^{3/4} \qquad (4.38)\\
&\qquad + 2\mathbb{I}_{\max_{0\le k\le n}\sup_{0\le t\le D\sqrt{n^{-1}\ln n}}|B_{\theta^{(n)}_k + t} - B_{\theta^{(n)}_k}| > Dn^{-1/4}(\ln n)^{3/4}} \\
&\qquad\quad \times \sup_{0\le t\le\theta^{(n)}_n} |B_t|.
\end{aligned}
$$

Since the sequence $\{\theta^{(n)}_l, l\ge 1\}$ has the same distribution as the sequence $\{\frac{T}{n}\Psi_l, l\ge 1\}$ defined in the proof of Lemma 4.1, then in the same way as in (4.16) and (4.17), we obtain

$$
\begin{aligned}
&P^B\left\{\max_{0\le l\le n}\left|\theta^{(n)}_l - \frac{lT}{n}\alpha_n\right| > \frac{1}{3}D\sqrt{n^{-1}\ln n} - K_4n^{-1}\right\} \qquad (4.39)\\
&\quad = P^B\left\{\max_{0\le l\le n}|\Psi_l - l\alpha_n| > \frac{D}{3T}\sqrt{n\ln n} - K_4T^{-1}\right\} \le e^{3K_4}n^{-1},
\end{aligned}
$$

provided we choose $D \ge 2 + 18T^2e^{2b_T}b_T^{-3}$ and assume that $n \ge 2(36T^2/b_T^2)^2$, so that $n^{-1}\ln n \le \frac{1}{36}T^{-2}b_T^2$. Hence, by (4.7), (4.9), (4.20), (4.37)–(4.39), the Cauchy–Schwarz and Chebyshev inequalities, it follows that there exists a constant $\tilde{C} > 0$ such that

$$P^B\{H_2^{(n)} > 1\} \le \tilde{C} n^{-1/2} (\ln n)^{3/2}.$$

This together with (4.7)–(4.9), (4.20), (4.25), (4.34), (4.37)–(4.39) and the Cauchy–Schwarz inequality yields that there exists a constant $C^{(2)} > 0$ such that

$$E^B J_{23} \le C^{(2)} z n^{-1/4} (\ln n)^{3/4} \tag{4.40}$$

and both $\tilde{C}$ and $C^{(2)}$ can be estimated explicitly from the above formulas. Finally, estimating the left-hand side of (4.28) by means of (4.7), (4.8), (4.25), (4.28), (4.29), (4.33), (4.40) and estimating the other two terms in the right-hand side of (4.27) exactly in the same way, we obtain that

$$J_2(\zeta, \eta) \le C^{(3)} z n^{-1/4} (\ln n)^{3/4} \tag{4.41}$$

for some $C^{(3)} > 0$ independent of $n$, which together with (4.22)–(4.26) yields (3.17), completing the proof of Lemma 3.2. □

Proof of Lemma 3.3. In order to prove Lemma 3.3, we write by (2.1), (2.9) and (3.15) that, for any $k, l = 1, 2, \ldots, n$,

$$\begin{aligned} &|tQ_z^{B,n,\theta}(\theta_k^{(n)}, \theta_l^{(n)}) - Q_z^B(\theta_k^{(n)}, \theta_l^{(n)})| \\ &\quad \le |R_z^{B,n,\theta}(\theta_k^{(n)}, \theta_l^{(n)}) - R_z^B(\theta_k^{(n)}, \theta_l^{(n)})| \\ &\quad \le L((\theta_n^{(n)} + 1) J_3), \end{aligned} \tag{4.42}$$

where

$$J_3 = \sup_{0 \le t \le \theta_n^{(n)}} |S_t^{B,\theta,n}(z) - S_t^B(z)|.$$

Set

$$H_3^{(n,l)} = \max_{1 \le k \le l} \sup_{\theta_{k-1}^{(n)} \le s \le t \le \theta_k^{(n)}} (r(t-s) + \kappa |B_t^* - B_s^*|) \quad \text{and} \quad H_3^{(n)} = H_3^{(n,n)}.$$

Then by (2.4) and (3.14), similarly to (4.19),

$$J_3 \le 2 z e^{r \theta_n^{(n)}} H_3^{(n)} \max_{0 \le k \le n} e^{\kappa B_{\theta_k^{(n)}}} + 2 z e^{r \theta_n^{(n)}} \mathbb{I}_{H_3^{(n)} > 1} \max_{0 \le k \le n} e^{\kappa B_{\theta_k^{(n)}}}. \tag{4.43}$$

Since $|B_t - B_s| \le |B_t - B_{\theta_{k-1}^{(n)}}| + |B_t - B_{\theta_{k-1}^{(n)}}|$ and $(a+b)^{2m} \le 2^{2m-1}(a^{2m} + b^{2m})$ for $a, b \ge 0$, $m \ge 1/2$, we obtain by (4.1), (4.5) and (4.9) that

$$E^B |H_3^{(n)}|^{2m} \le 2^{2m-1} \left| r - \frac{\kappa^2}{2} \right|^{2m} E^B \max_{1 \le k \le n} |\theta_k^{(n)} - \theta_{k-1}^{(n)}|^{2m}$$

$$
\begin{aligned}
&+ 2^{4m-1}\kappa^{2m} E^B \max_{1\le k\le n} \sup_{\theta_{k-1}^{(n)}\le t\le \theta_k^{(n)}} |B_t - B_{\theta_{k-1}^{(n)}}|^{2m} \\
&\le 2^{2m-1}\sum_{k=1}^{n}\biggl(\left|r-\frac{\kappa^2}{2}\right|^{2m} E^B|\theta_k^{(n)}-\theta_{k-1}^{(n)}|^{2m} \\
&\qquad + 2^{2m}\kappa^{2m}E^B \sup_{\theta_{k-1}^{(n)}\le t\le \theta_k^{(n)}} |B_t - B_{\theta_{k-1}^{(n)}}|^{2m}\biggr) \\
&\le 2^{2m-1} n\left(\left|r-\frac{\kappa^2}{2}\right|^{2m} E^B(\theta_1^{(n)})^{2m} + 2^{2m}\kappa^{2m}\Lambda_m E^B(\theta_1^{(n)})^m\right) \\
&\le K_5^{(m)} n^{-m+1},
\end{aligned}
\tag{4.44}
$$

where

$$
K_5^{(m)} = 2^{2m-1}e^{2b_T} b_T^{-(m+1)} T^m \left(\left|r-\frac{\kappa^2}{2}\right|^{2m}(2m)!b_T^{-m}T^m + 2^{2m}\kappa^{2m}\Lambda_m m!\right).
$$

This together with the Chebyshev inequality gives that, for any integers $m, n \ge 1$,

$$
P^B\{H_3^{(n)} > 1\} \le K_5^{(m)} n^{-m+1}. \tag{4.45}
$$

Finally, by (4.8), (4.25), (4.42)–(4.45) and the Hölder inequality, we obtain the assertion of Lemma 3.3. □

Proof of Lemma 3.4. The proof of Lemma 3.4 starts similarly with the estimate

$$
|Q_z^B(\theta_k^{(n)}, \theta_l^{(n)}) - Q_z^B(\theta_k^{(n)}\wedge T, \theta_l^{(n)}\wedge T)| \le E^B(J_4(k,l) + J_5(k,l)), \tag{4.46}
$$

where by (2.2)–(2.4) and (2.8),

$$
\begin{aligned}
J_4(k,l) &= |e^{-r\theta_{k\wedge l}^{(n)}} - e^{-r\theta_{k\wedge l}^{(n)}\wedge T}|R_z^B(\theta_k^{(n)}\wedge T, \theta_l^{(n)}\wedge T) \\
&\le rL(T+2)(F_0(z)+\Delta_0(z)+z+1)e^{rT}\left(|\theta_n^{(n)}-T| \sup_{0\le t\le T} e^{\kappa B_t}\right)
\end{aligned}
\tag{4.47}
$$

and

$$
\begin{aligned}
J_5(k,l) &\le |R_z^B(\theta_k^{(n)}, \theta_l^{(n)}) - R_z^B(\theta_k^{(n)}\wedge T, \theta_l^{(n)}\wedge T)| \\
&\le L\left(|\theta_n^{(n)}-T|\left(1+e^{r\theta_n^{(n)}} \sup_{0\le t\le \theta_n^{(n)}} e^{\kappa B_t}\right) + J_{51}\right)
\end{aligned}
\tag{4.48}
$$

with

$$
J_{51} = \max_{0\le k\le n} \sup_{\theta_k^{(n)}\wedge T\le t\le \theta_k^{(n)}} |S_t^B(z) - S_{\theta_k^{(n)}\wedge T}^B(z)|.
$$

Set

$$H_4^{(n)} = \max_{0\le k\le n} \sup_{\theta_k^{(n)}\wedge T\le t\le \theta_k^{(n)}} (r(t-\theta_k^{(n)}\wedge T)+\kappa|B_t^*-B^*_{\theta_k^{(n)}\wedge T}|).$$

By (2.4), similarly to (4.30) and (4.43), we obtain that

$$\begin{aligned} J_{51} &\le ze^{rT}H_4^{(n)} \sup_{0\le t\le T} e^{\kappa B_t} \\ &\quad + z\mathbb{I}_{H_4^{(n)}>1}\Big(e^{r\theta_n^{(n)}} \sup_{0\le t\le \theta_n^{(n)}} e^{\kappa B_t} + e^{rT} \sup_{0\le t\le T} e^{\kappa B_t}\Big). \end{aligned} \tag{4.49}$$

Observe that $\theta_k^{(n)}\wedge T<\theta_k^{(n)}$ if and only if $T<\theta_k^{(n)}$ and since $k\le n$, we have in this case

$$[\theta_k^{(n)}\wedge T, \theta_k^{(n)}]\subset [T,\theta_n^{(n)}]\subset [T,\theta_n^{(n)}\vee T].$$

Hence,

$$H_4^{(n)} \le \Big(r+\frac{\kappa^2}{2}\Big)|\theta_n^{(n)}-T| + \kappa \sup_{T\le t\le \theta_n^{(n)}\vee T} |B_t-B_T|. \tag{4.50}$$

By (4.7), (4.9) and the Cauchy–Schwarz inequality,

$$E^B \sup_{T\le t\le \theta_n^{(n)}\vee T} |B_t-B_T|^2 \le \Lambda_1 E^B|\theta_n^{(n)}-T| \le \Lambda_1 T K_2^{1/2} n^{-1/2},$$

which together with (4.7) and (4.50) gives

$$E^B|H_4^{(n)}|^2 \le 2\Big(r+\frac{\kappa^2}{2}\Big)^2 K^2T^2n^{-1} + 2\kappa^2\Lambda_1 T K_2^{1/2} n^{-1/2}.$$

This enables us to estimate $P^B\{H_4^{(n)}>1\}$ by the Chebyshev inequality and together with (4.8), (4.25) and the Cauchy–Schwarz inequality yields

$$E^B J_{51} \le C^{(4)} z n^{-1/4}$$

for some $C^{(4)}>0$ independent of $n$ which can be easily estimated explicitly via the above formulas. Using (4.7), (4.8) and (4.25) together with the Cauchy–Schwarz inequality in order to estimate $E^B J_4(\zeta,\eta)$ and the expectation of the remaining term in $J_5(\zeta,\eta)$, we arrive at (3.19), completing the proof of Lemma 3.4. □

Next, we derive Lemma 3.5 using estimates similar to the above.

Proof of Lemma 3.5. Namely, for any $k,l=1,2,\dots,$ we have

$$\begin{aligned} &|R_z^B(\theta_k^{(n)}\wedge T, \theta_l^{(n)}\wedge T) - R_z^B(\theta_{k\wedge n}^{(n)}\wedge T, \theta_{l\wedge n}^{(n)}\wedge T)| \\ &\quad \le E^B(J_6(k,l)+J_7(k,l)), \end{aligned} \tag{4.51}$$

where by (2.2)–(2.4), (2.9) and the equality $\theta_\zeta^{(n)} \wedge \theta_\eta^{(n)} = \theta_{\zeta\wedge\eta}^{(n)}$ similarly to (4.43) and (4.49),

$$\begin{aligned} J_6(k,l) &= (|e^{-r\theta_{k\wedge l}^{(n)}\wedge T} - e^{-r\theta_{k\wedge l\wedge n}^{(n)}\wedge T}|R_z^B(\theta_k^{(n)} \wedge T, \theta_l^{(n)} \wedge T)) \\ &\le rL(T+2)(F_0(z)+\Delta_0(z)+z+1)e^{rT} \\ &\quad \times \max_{n<k<\infty} |\theta_k^{(n)} \wedge T - \theta_n^{(n)} \wedge T| \sup_{0\le t\le T} e^{\kappa B_t} \end{aligned} \tag{4.52}$$

and

$$\begin{aligned} J_7(k,l) &\le |R_z^B(\theta_k^{(n)} \wedge T, \theta_l^{(n)} \wedge T) - R_z^B(\theta_{k\wedge n}^{(n)} \wedge T, \theta_{l\wedge n}^{(n)} \wedge T)| \\ &\le Lze^{rT}\biggl(\Bigl(1 + \sup_{0\le t\le T} e^{\kappa B_t}\Bigr) \\ &\qquad \times \Bigl(\max_{n<k<\infty} |\theta_k^{(n)} \wedge T - \theta_n^{(n)} \wedge T| + H_5^{(n)} + 2\mathbb{I}_{H_5^{(n)}>1}\Bigr)\biggr), \end{aligned} \tag{4.53}$$

with

$$H_5^{(n)} = \max_{n<k<\infty} \sup_{\theta_n^{(n)}\wedge T\le t\le \theta_k^{(n)}\wedge T} (r(t-\theta_n^{(n)}\wedge T) + \kappa|B_t^* - B^*_{\theta_n^{(n)}\wedge T}|).$$

Observe that $\theta_n^{(n)} \wedge T < \theta_k^{(n)} \wedge T$ for $k > n$ if and only if $T > \theta_n^{(n)}$ and then

$$[\theta_n^{(n)} \wedge T, \theta_k^{(n)} \wedge T] \subset [\theta_n^{(n)}, T] \subset [\theta_n^{(n)}, \theta_n^{(n)} \vee T].$$

Hence,

$$\begin{aligned} |\theta_k^{(n)} \wedge T - \theta_n^{(n)} \wedge T| &\le T \vee \theta_n^{(n)} - \theta_n^{(n)} \\ &\le |T - \theta_n^{(n)}| \end{aligned} \tag{4.54}$$

and

$$H_5^{(n)} \le \Bigl(r + \frac{\kappa^2}{2}\Bigr)|T - \theta_n^{(n)}| + \kappa \sup_{T\le t\le \theta_n^{(n)}\vee T} |B_t - B_T|. \tag{4.55}$$

The right-hand side of (4.55) is the same as in (4.50), and so we can use the same estimates for $H_5^{(n)}$ as for $H_4^{(n)}$, which together with (4.7), (4.8), (4.25), the Chebyshev and the Cauchy–Schwarz inequalities enable us to estimate $E^B J_7(k,l)$ by $C^{(5)} z n^{-1/4}$ for some $C^{(5)} > 0$ independent of $n$. Finally, using (4.7), (4.8), (4.25) and (4.54) together with the Cauchy–Schwarz inequality in order to estimate $E^B J_6(k,l)$, we obtain (3.20), completing the proof of Lemma 3.5. □

In order to complete the proof of Theorem 2.1, it remains to establish Lemma 3.6.

PROOF OF LEMMA 3.6. For each $\sigma \in \mathcal{T}^B_{0,T}$, set $\nu_\sigma = \min\{k \in \mathbb{N} : \theta^{(n)}_k \geq \sigma\}$ which, indeed, defines $\nu_\sigma$ since $\theta^{(n)}_k \to \infty$ with probability one as $k \to \infty$. Observe that $\nu_\sigma \in \mathcal{T}^{B,n}$ since $\{\nu_\sigma \leq k\} = \{\theta^{(n)}_k \geq \sigma\} \in \mathcal{F}^B_{\theta^{(n)}_k}$. For any $\sigma \in \mathcal{T}^B_{0,T}$, we set $\sigma^{(n)} = \theta^{(n)}_{\nu_\sigma} \wedge T$. Since $\mathcal{T}^{B,n}_T \subset \mathcal{T}^B_{0,T}$, we conclude from (2.10) that

$$V(z) \geq \inf_{\sigma \in \mathcal{T}^B_{0,T}} \sup_{\tau \in \mathcal{T}^{B,n}_T} E^B Q^B_z(\sigma, \tau). \tag{4.56}$$

Then for any $\delta > 0$, there exists $\sigma_\delta \in \mathcal{T}^B_{0,T}$ such that

$$V(z) \geq \sup_{\tau \in \mathcal{T}^{B,n}_T} E^B Q^B_z(\sigma_\delta, \tau) - \delta. \tag{4.57}$$

This together with (3.21) implies

$$\begin{aligned} V(z) &\geq \sup_{\tau \in \mathcal{T}^{B,n}_T} E^B Q^B_z(\sigma^{(n)}_\delta, \tau) - \delta \\ &\quad - \sup_{\tau \in \mathcal{T}^{B,n}_T} E^B(Q^B_z(\sigma^{(n)}_\delta, \tau) - Q^B_z(\sigma_\delta, \tau)) \\ &\geq V^{B,n}_{0,T} - \delta - \sup_{\tau \in \mathcal{T}^{B,n}_T} J_8(\sigma_\delta, \tau) - \sup_{\tau \in \mathcal{T}^{B,n}_T} J_9(\sigma_\delta, \tau), \end{aligned} \tag{4.58}$$

where, for any $\sigma \in \mathcal{T}^B_{0,T}$ and $\tau \in \mathcal{T}^{B,n}_T$,

$$J_8(\sigma, \tau) = E^B(e^{-r\sigma^{(n)} \wedge \tau}(R^B_z(\sigma^{(n)}, \tau) - R^B_z(\sigma, \tau)))$$

and by (2.3)–(2.4),

$$\begin{aligned} J_9(\sigma, \tau) &= E^B(|e^{-r\sigma^{(n)} \wedge \tau} - e^{-r\sigma \wedge \tau}| R^B_z(\sigma, \tau)) \\ &\leq rL(T+2)(F_0(z) + \Delta_0(z) + z)e^{rT} \\ &\quad \times E^B\Big(\sup_{0 \leq k < \infty} |\theta^{(n)}_{k+1} \wedge T - \theta^{(n)}_k \wedge T| \Big(1 + \sup_{0 \leq t \leq T} e^{\kappa B_t}\Big)\Big). \end{aligned} \tag{4.59}$$

Since $\sigma^{(n)} \geq \sigma$, it follows that

$$R^B_z(\sigma, \tau) = F_\sigma(S^B(z)) + \Delta_\sigma(S^B(z))$$

whenever

$$R^B_z(\sigma^{(n)}, \tau) = F_{\sigma^{(n)}}(S^B(z)) + \Delta_{\sigma^{(n)}}(S^B(z)).$$

Thus, by (2.2) and (2.8), similarly to (4.53) for any $\sigma \in \mathcal{T}^B_{0T}$ and $\tau \in \mathcal{T}^{B,n}_T$,

$$R^B_z(\sigma^{(n)}, \tau) - R^B_z(\sigma, \tau)$$

$$
\begin{aligned}
&\le |F_{\sigma^{(n)}}(S^B(z)) - F_\sigma(S^B(z))| + |\Delta_{\sigma^{(n)}}(S^B(z)) - \Delta_\sigma(S^B(z))| \\
&\le Lze^{rT}\left(1+\sup_{0\le t\le T} e^{\kappa B_t}\right) \\
&\quad \times \left(\sup_{0\le k<\infty}|\theta_{k+1}^{(n)}\wedge T - \theta_k^{(n)}\wedge T| + H_6^{(n)} + 2\mathbb{I}_{H_6^{(n)}>1}\right),
\end{aligned}
\tag{4.60}
$$

where

$$
H_6^{(n)} = \sup_{0\le k<\infty} \sup_{\theta_k^{(n)}\wedge T\le t\le \theta_{k+1}^{(n)}\wedge T} (r(t-\theta_k^{(n)}\wedge T) + \kappa(B_t^* - B^*_{\theta_k^{(n)}\wedge T})).
$$

It is clear that

$$
|H_6^{(n)}| \le |H_3^{(n)}| + |H_5^{(n)}| \tag{4.61}
$$

and by (4.54), for all $k\ge 0$,

$$
|\theta_{k+1}^{(n)}\wedge T - \theta_k^{(n)}\wedge T| \le \max_{0\le k\le n-1}|\theta_{k+1}^{(n)} - \theta_k^{(n)}| + |T-\theta_n^{(n)}|. \tag{4.62}
$$

Hence, we can apply to the right-hand side of (4.60) the estimates of Lemmas 3.3 and 3.5 arriving at a bound of order $n^{-1/4}$. In order to obtain a better estimate promised in Lemma 3.6 (though it will not help us to improve the estimate of Theorem 2.1), we write

$$
|H_6^{(n)}| \le |H_3^{(n,2n)}| + \mathbb{I}_{\theta_{2n}^{(n)}<T}(|H_3^{(n)}| + |H_5^{(n)}|), \tag{4.63}
$$

with $H_3^{(n,l)}$ defined above (4.43). In the same way as in (4.44), we obtain

$$
E^B|H_3^{(n,2n)}|^{2m} \le 2K_5^{(m)} n^{-m+1}. \tag{4.64}
$$

Next, by (4.1) and (4.4), similarly to (4.17), (4.18) and (4.39), we obtain the following (large deviations) upper bound:

$$
P^B\{\theta_{2n}^{(n)} < T\} = P^B\{2n\alpha_n - \Psi_{2n} > n(2\alpha_n - 1)\} \le \tilde{C}e^{-\rho n} \tag{4.65}
$$

for all $n\in\mathbb{N}$ and some $\tilde{C},\rho>0$ independent of $n$. Estimating $E^B|T-\theta_n^{(n)}|^m$ by (4.7), we obtain by (4.9), (4.44), (4.55), (4.63)–(4.65) and the Cauchy–Schwarz inequality that, for any $m\ge 1$, there exists $K_6^{(m)}>0$ (which can be explicitly estimated from above formulas) such that

$$
E^B|H_6^{(n)}|^{2m} \le K_6^{(m)} n^{-m+1}. \tag{4.66}
$$

Since $\delta$ in (4.58) is arbitrary, we conclude by (4.7), (4.8), (4.25), (4.55), (4.58)–(4.60), (4.62) and (4.66) together with the Chebyshev and Hölder inequalities similarly to Lemma 3.3 that, for any $\varepsilon>0$, there exists $C_\varepsilon^{(6)}>0$ such that, for all $n\in\mathbb{N}$,

$$
V(z) - V_{0,T}^{B,n}(z) \ge -C_\varepsilon^{(6)}(F_0(z) + \Delta_0(z) + z + 1)n^{\varepsilon-1/2}. \tag{4.67}
$$

By [23], we can represent $V(z)$ not only in the form (2.10), but also as

$$
\begin{aligned}
V(z) &= \sup_{\tau\in\mathcal{T}_{0T}^B}\inf_{\sigma\in\mathcal{T}_{0T}^B} E^B Q_z^B(\sigma,\tau)\\
&\le \inf_{\sigma\in\mathcal{T}_T^{B,n}} E^B Q_z^B(\sigma,\tau_\delta)+\delta
\end{aligned}
$$

for each $\delta>0$ and some $\tau_\delta\in\mathcal{T}_{0T}^B$. Introducing $\tau_\delta^{(n)}$ and employing the same arguments as above, we obtain that, for any $\varepsilon>0$, there exists $C_\varepsilon^{(7)}>0$ such that, for all $n\in\mathbb{N}$,

$$
V(z)-V_{0,T}^{B,n}(z)\le C_\varepsilon^{(7)}(F_0(z)+\Delta_0(z)+z+1)n^{\varepsilon-1/2},
$$

which together with (4.67) yields (3.23) and completes the proof of both Lemma 3.6 and Theorem 2.1. $\square$

**5. Exercise times and hedging with small shortfalls.** Set

$$
\begin{aligned}
\underline{V}^{(n)}(z,\eta) &= \min_{\zeta\in\mathcal{T}_{0n}^\xi} E_n^\xi Q_z^{(n)}\left(\frac{\zeta T}{n},\frac{\eta T}{n}\right),\\
\overline{V}^{(n)}(z,\zeta) &= \max_{\eta\in\mathcal{T}_{0n}^\xi} E_n^\xi Q_z^{(n)}\left(\frac{\zeta T}{n},\frac{\eta T}{n}\right),
\end{aligned}
\tag{5.1}
$$

$$
\begin{aligned}
\underline{V}^{B,n}(z,\eta) &= \inf_{\zeta\in\mathcal{T}_{0,n}^{B,n}} E^B Q_z^{B,n}\left(\frac{\zeta T}{n},\frac{\eta T}{n}\right),\\
\overline{V}^{B,n}(z,\zeta) &= \sup_{\eta\in\mathcal{T}_{0,n}^{B,n}} E^B Q_z^{B,n}\left(\frac{\zeta T}{n},\frac{\eta T}{n}\right)
\end{aligned}
\tag{5.2}
$$

and

$$
\begin{aligned}
\underline{V}_{\mathcal{S}}^{B,n}(z,\eta) &= \min_{\zeta\in\mathcal{S}_{0,n}^{B,n}} E^B Q_z^{B,n}\left(\frac{\zeta T}{n},\frac{\eta T}{n}\right),\\
\overline{V}_{\mathcal{S}}^{B,n}(z,\zeta) &= \max_{\eta\in\mathcal{S}_{0,n}^{B,n}} E^B Q_z^{B,n}\left(\frac{\zeta T}{n},\frac{\eta T}{n}\right).
\end{aligned}
\tag{5.3}
$$

We define similar quantities $\hat{\underline{V}}^{(n)}, \hat{\overline{V}}^{(n)}, \hat{\underline{V}}^{B,n}, \hat{\overline{V}}^{B,n}, \hat{\underline{V}}_{\mathcal{S}}^{B,n}, \hat{\overline{V}}_{\mathcal{S}}^{B,n}$ replacing $Q_z^{(n)}$ by $\hat{Q}_z^{(n)}$ and $Q_z^{B,n}$ by $\hat{Q}_z^{B,n}$ in the corresponding formulas.

Proof of Theorem 2.2. The starting point in the proof of Theorem 2.2 is the following result similar to Lemma 3.1.

LEMMA 5.1. *For any $z, n > 0$ and each $\mu \in \mathcal{J}_{0n}$,*

$$\underline{V}_{\mathcal{S}}^{B,n}(z, \mu \circ \lambda_B^{(n)}) = \underline{V}^{(n)}(z, \mu \circ \lambda_\xi^{(n)}) = \underline{V}^{B,n}(z, \mu \circ \lambda_B^{(n)}) \tag{5.4}$$

*and*

$$\overline{V}_{\mathcal{S}}^{B,n}(z, \mu \circ \lambda_B^{(n)}) = \overline{V}^{(n)}(z, \mu \circ \lambda_\xi^{(n)}) = \overline{V}^{B,n}(z, \mu \circ \lambda_B^{(n)}). \tag{5.5}$$

*The corresponding results hold true also for $\hat{\underline{V}}^{(n)}, \hat{\overline{V}}^{(n)}, \hat{\underline{V}}^{B,n}, \hat{\overline{V}}^{B,n}, \hat{\underline{V}}_{\mathcal{S}}^{B,n}, \hat{\overline{V}}_{\mathcal{S}}^{B,n}$ in place of $\underline{V}^{(n)}, \overline{V}^{(n)}, \underline{V}^{B,n}, \overline{V}^{B,n}, \underline{V}_{\mathcal{S}}^{B,n}, \overline{V}_{\mathcal{S}}^{B,n}$, respectively.*

PROOF. The first equality in (5.4) and (5.5) follows in the same way as the first equality in (3.5), taking into account that

$$Q_z^{(n)}(\mu \circ \lambda_\xi^{(n)}, \nu \circ \lambda_\xi^{(n)}) \quad \text{and} \quad Q_z^{B,n}(\mu \circ \lambda_B^{(n)}, \nu \circ \lambda_B^{(n)})$$

have the same distribution for all $\mu, \nu \in \mathcal{J}_{0n}$. In order to obtain the second equality in (5.4) and (5.5), we employ again the well-known dynamical programming recursive relations for the optimal stopping problem (see, e.g., [26]) which have here the form

$$\begin{aligned} &\underline{V}^{(n)}(z, \eta) = \underline{V}_{0,n}^{(n)}(z, \eta), \qquad \underline{V}_{n,n}^{(n)}(z, \eta) = Q_z^{(n)}\left(T, \frac{\eta T}{n}\right) \quad \text{and} \\ &\underline{V}_{k,n}^{(n)}(z, \eta) = \min\left(Q_z^{(n)}\left(\frac{kT}{n}, \frac{\eta T}{n}\right), E^\xi(V_{k+1,n}^{(n)}(z, \eta) | \mathcal{F}_k^\xi)\right), \end{aligned} \tag{5.6}$$

$$\begin{aligned} &\overline{V}^{(n)}(z, \zeta) = \overline{V}_{0,n}^{(n)}(z, \zeta), \qquad \overline{V}_{n,n}^{(n)}(z, \zeta) = Q_z^{(n)}\left(\frac{\zeta T}{n}, T\right) \quad \text{and} \\ &\overline{V}_{k,n}^{(n)}(z, \eta) = \max\left(Q_z^{(n)}\left(\frac{\zeta T}{n}, \frac{kT}{n}\right), E^\xi(V_{k+1,n}^{(n)}(z, \zeta) | \mathcal{F}_k^\xi)\right) \end{aligned} \tag{5.7}$$

for any $\zeta, \eta \in \mathcal{T}_{0n}^\xi$ and

$$\begin{aligned} &\underline{V}^{B,n}(z, \eta) = \underline{V}_{0,n}^{B,n}(z, \eta), \qquad \underline{V}_{n,n}^{B,n}(z, \eta) = Q_z^{B,n}\left(T, \frac{\eta T}{n}\right) \quad \text{and} \\ &\underline{V}_{k,n}^{B,n}(z, \eta) = \min\left(Q_z^{B,n}\left(kTn, \frac{\eta T}{n}\right), E^B(V_{k+1,n}^{B,n}(z, \eta) | \mathcal{F}_{\theta_k^{(n)}}^B)\right), \end{aligned} \tag{5.8}$$

$$\begin{aligned} &\overline{V}^{B,n}(z, \zeta) = \overline{V}_{0,n}^{B,n}(z, \zeta), \qquad \overline{V}_{n,n}^{B,n}(z, \zeta) = Q_z^{B,n}\left(\frac{\zeta T}{n}, T\right) \quad \text{and} \\ &\overline{V}_{k,n}^{B,n}(z, \eta) = \max\left(Q_z^{B,n}\left(\frac{\zeta T}{n}, \frac{kT}{n}\right), E^B(V_{k+1,n}^{B,n}(z, \zeta) | \mathcal{F}_{\theta_k^{(n)}}^B)\right) \end{aligned} \tag{5.9}$$

for any $\zeta, \eta \in \mathcal{T}_{0n}^\xi$.

It is clear from the construction of the stopping times $\mu\circ\lambda_\xi^{(n)}$ and $\mu\circ\lambda_B^{(n)}$ for $\mu\in J_{0n}$ that $(\mu\circ\lambda_\xi^{(n)})\wedge k$, $\mathbb{I}_{\mu\circ\lambda_\xi^{(n)}\geq k}$ and $(\mu\circ\lambda_B^{(n)})\wedge k$, $\mathbb{I}_{\mu\circ\lambda_B^{(n)}\geq k}$ are measurable with respect to the $\sigma$-algebras $\mathcal{F}_k^\xi$ and $\mathcal{G}_k^{B,n}$, respectively. Since $F_k(S^{(n)})$, $\Delta_k(S^{(n)})$ and $F_k(S^B)$, $\Delta_k(S^B)$ are also $\mathcal{F}_k^\xi$- and $\mathcal{G}_k^{B,n}$-measurable, respectively, we conclude that $Q_z^{(n)}(\frac{kT}{n},\frac{T}{n}\mu\circ\lambda_\xi^{(n)})$, $Q_z^{(n)}(\frac{T}{n}\mu\circ\lambda_\xi^{(n)},\frac{kT}{n})$ and $Q_z^{B,n}(\frac{kT}{n},\frac{T}{n}\mu\circ\lambda_B^{(n)})$, $Q_z^{B,n}(\frac{T}{n}\mu\circ\lambda_B^{(n)},\frac{kT}{n})$ are $\mathcal{F}_k^\xi$- and $\mathcal{G}_k^{B,n}$-measurable, respectively. It follows from here by the backward induction in the same way as in Lemma 3.1 that there exist measurable functions $\underline{\Phi}_k^\mu(z,x_1,\dots,x_k)$, $\underline{\Phi}_0^\mu(z)$ and $\overline{\Phi}_k^\mu(z,x_1,\dots,x_k)$, $\overline{\Phi}_0^\mu(z)$, $k=1,2,\dots,n$, such that

$$\begin{aligned}&\underline{V}_{kn}^{(n)}(z,\mu\circ\lambda_\xi^{(k)})=\underline{\Phi}_k^\mu\Big(z,\Big(\frac{T}{n}\Big)^{1/2}\xi_1,\dots,\Big(\frac{T}{n}\Big)^{1/2}\xi_k\Big),\\&\underline{V}_{0n}^{(n)}(z,\mu\circ\lambda_\xi^{(n)})=\underline{\Phi}_0^\mu(z),\end{aligned}\tag{5.10}$$

$$\begin{aligned}&\overline{V}_{kn}^{(n)}(z,\mu\circ\lambda_\xi^{(k)})=\overline{\Phi}_k^\mu\Big(z,\Big(\frac{T}{n}\Big)^{1/2}\xi_1,\dots,\Big(\frac{T}{n}\Big)^{1/2}\xi_k\Big)\\&\overline{V}_{0n}^{(n)}(z,\mu\circ\lambda_\xi^{(n)})=\overline{\Phi}_0^\mu(z)\end{aligned}\tag{5.11}$$

and

$$\begin{aligned}&\underline{V}_{kn}^{B,n}(z,\mu\circ\lambda_B^{(k)})=\underline{\Phi}_k^\mu(z,B^*_{\theta_1^{(n)}},B^*_{\theta_2^{(n)}}-B^*_{\theta_1^{(n)}},\dots,B^*_{\theta_k^{(n)}}-B^*_{\theta_{k-1}^{(n)}}),\\&\underline{V}_{0n}^{B,n}(z,\mu\circ\lambda_B^{(n)})=\underline{\Phi}_0^\mu(z),\end{aligned}\tag{5.12}$$

$$\begin{aligned}&\overline{V}_{kn}^{B,n}(z,\mu\circ\lambda_B^{(k)})=\overline{\Phi}_k^\mu(z,B^*_{\theta_1^{(n)}},B^*_{\theta_2^{(n)}}-B^*_{\theta_1^{(n)}},\dots,B^*_{\theta_k^{(n)}}-B^*_{\theta_{k-1}^{(n)}}),\\&\overline{V}_{0n}^{B,n}(z,\mu\circ\lambda_B^{(n)})=\overline{\Phi}_0^\mu(z).\end{aligned}\tag{5.13}$$

Applying these formulas for $k=0$, we obtain the second equality in (5.4) and (5.5). The corresponding results for $\hat{\underline{V}}^{(n)},\hat{\overline{V}}^{(n)},\hat{\underline{V}}^{B,n},\hat{\overline{V}}^{B,n},\hat{\underline{V}}_{\mathcal{S}}^{B,n},\hat{\overline{V}}_{\mathcal{S}}^{B,n}$ are derived in the same way. □

Define

$$\begin{aligned}&\underline{V}_{0,T}^{B,n}(z,\tau)=\inf_{\sigma\in\mathcal{T}_T^{B,n}}E^BQ_z^B(\sigma,\tau)\quad\text{and}\\&\overline{V}_{0,T}^{B,n}(z,\sigma)=\sup_{\tau\in\mathcal{T}_T^{B,n}}E^BQ_z^B(\sigma,\tau).\end{aligned}\tag{5.14}$$

The proof of Lemma 3.2, together with Lemmas 3.3–3.5, yields that there exists $C^{(8)}>0$ such that, for each $\mu\in\mathcal{J}_{0n}$ and all $z,n>0$,

$$\begin{aligned}&|\underline{V}^{B,n}(z,\mu\circ\lambda_B^{(n)})-\underline{V}_{0,T}^{B,n}(z,\theta^{(n)}_{\mu\circ\lambda_B^{(n)}}\wedge T)|\\&\quad\le C^{(8)}(F_0(z)+\Delta_0(z)+z+1)n^{-1/4}(\ln n)^{3/4}\end{aligned}\tag{5.15}$$

and

$$\begin{aligned}&|\overline{V}^{B,n}(z,\mu\circ\lambda_B^{(n)})-\overline{V}_{0,T}^{B,n}(z,\theta^{(n)}_{\mu\circ\lambda_B^{(n)}}\wedge T)|\\&\quad\le C^{(8)}(F_0(z)+\Delta_0(z)+z+1)n^{-1/4}(\ln n)^{3/4}.\end{aligned}\tag{5.16}$$

Next, set

$$\underline{V}(z,\tau)=\inf_{\sigma\in\mathcal{T}_{0T}^B}E^BQ_z^B(\sigma,\tau),\qquad \overline{V}(z,\sigma)=\sup_{\tau\in\mathcal{T}_{0T}^B}E^BQ_z^B(\sigma,\tau),\tag{5.17}$$

and similarly to Lemma 3.6, estimate $|\underline{V}_{0,T}^{B,n}(z,\tau)-\underline{V}(z,\tau)|$ and $|\overline{V}_{0,T}^{B,n}(z,\tau)-\overline{V}(z,\tau)|$, but now it is a bit simpler since we have obvious one-sided inequalities

$$\underline{V}_{0,T}^{B,n}(z,\tau)\ge\underline{V}(z,\tau)\quad\text{and}\quad\overline{V}_{0,T}^{B,n}(z,\tau)\le\overline{V}(z,\tau).\tag{5.18}$$

For any $\tau\in\mathcal{T}_T^{B,n}$ and $\delta>0$, there exists $\sigma_\delta\in\mathcal{T}_{0,T}^B$ such that

$$\underline{V}(z,\tau)=E^BQ_z^B(\sigma_\delta,\tau)-\delta,$$

and so using the notation $J_8$ and $J_9$ from (4.58) together with their estimates, we see that, for any $\varepsilon>0$, there exists $C_\varepsilon^{(9)}>0$ such that, for any $\delta>0$ and $n\in\mathbb{N}$,

$$\begin{aligned}\underline{V}(z,\tau)&=E^BQ_z^B(\sigma_\delta^{(n)},\tau)-\delta-J_8(\sigma_\delta,\tau)-J_9(\sigma_\delta,\tau)\\&\ge\underline{V}_{0,T}^{B,n}(z,\tau)-\delta-C_\varepsilon^{(9)}n^{\varepsilon-1/2}.\end{aligned}$$

Since $\delta$ is arbitrary and we have already inequality (5.18) in the other direction, it follows that

$$|\underline{V}(z,\tau)-\underline{V}_{0,T}^{B,n}(z,\tau)|\le C_\varepsilon^{(9)}n^{\varepsilon-1/2}.\tag{5.19}$$

Similarly, we obtain that

$$|\overline{V}(z,\tau)-\overline{V}_{0,T}^{B,n}(z,\tau)|\le C_\varepsilon^{(9)}n^{\varepsilon-1/2}.\tag{5.20}$$

Finally, (5.4), (5.5), (5.10)–(5.13), (5.19) and (5.20) together with (2.13) yield (2.17), provided (2.16) holds true. The proof is similar for $\hat\lambda_\xi^{(n)}$, $\hat\lambda_B^{(n)}$, $\hat Q_z^{(n)}$, $\hat V^{(n)}(z)$ in place of $\lambda_\xi^{(n)}$, $\lambda_B^{(n)}$, $Q_z^{(n)}$, $V^{(n)}(z)$. $\square$

Next, we establish Theorem 2.3.

Proof of Theorem 2.3. Since $\beta_k = \beta_k^\zeta$ and $\gamma_k = \gamma_k^\zeta$ in (2.19) are $\mathcal{F}_{k-1}$-measurable, they can be written uniquely in the form $\beta_k = f_k \circ \lambda_\xi^{(l)}$ and $\gamma_k = g_k \circ \lambda_\xi^{(l)}$ for any $l = k-1, k, \dots, n$, where $f_k$ and $g_k$ are considered as functions on $\{-1,1\}^n$ depending only on first $k-1$ variables, that is, in fact, as functions on $\{-1,1\}^{k-1}$.

In order to show that the portfolio $Z_t^B$ defined by (2.21) is self-financing, it suffices to check that the discounted portfolio

$$\check{Z}_t^B = Z_t^B b_t^{-1} = \beta_t^\varphi b_0 + \gamma_t^\varphi \check{S}_t^B, \qquad \check{S}_t^B = S_t^B e^{-rt} \tag{5.21}$$

is self-financing, which means in our case that, with probability one,

$$(f_{k+1} - f_k) \circ \lambda_B^{(k)} b_0 + (g_{k+1} - g_k) \circ \lambda_B^{(k)} \check{S}^B_{\theta_k^{(n)}} = 0, \tag{5.22}$$

where, recall, $\beta_t^\varphi = f_k \circ \lambda_B^{(k-1)} = f_k \circ \lambda_B^{(k)}$ and $\gamma_t^\varphi = g_k \circ \lambda_B^{(k-1)} = g_k \circ \lambda_B^{(k)}$ whenever $t \in (\theta_{k-1}^{(n)}, \theta_k^{(n)}]$. But if $\lambda_B^{(k)}(\omega) = \lambda_\xi^{(k)}(\omega')$, then

$$\begin{aligned} &(f_{k+1} - f_k) \circ \lambda_B^{(k)}(\omega) b_0 + (g_{k+1} - g_k) \circ \lambda_B^{(k)}(\omega) \check{S}^B_{\theta_k^{(n)}(\omega)}(z, \omega) \\ &\qquad = (f_{k+1} - f_k) \circ \lambda_\xi^{(k)}(\omega') b_0 + (g_{k+1} - g_k) \circ \lambda_\xi^{(k)}(\omega') \check{S}^{(n)}_{kT/n}(z, \omega'), \end{aligned} \tag{5.23}$$

where $\check{S}^{(n)}_{kT/n}(z) = S^{(n)}_{kT/n}(z) e^{-rkT/n}$. By (2.19), the right-hand side in (5.23) equals zero, and so (5.22) follows.

Recall that the sequence $B^*_{\theta_k^{(n)}} - B^*_{\theta_{k-1}^{(n)}}$, $k = 1, 2, \dots$, has the same distribution as the sequence $(\frac{T}{n})^{1/2} \xi_k$, $k = 1, 2, \dots$. Since the processes $S_t^{B,n}(z)$, $t \geq 0$, and $S_t^{(n)}(z)$, $t \geq 0$, defined by (3.1) and (2.5), respectively, are obtained by the same formulas from these sequences, they have the same distribution as well. Moreover, if $\mu \in \mathcal{J}_{0n}$, $\zeta = \mu \circ \lambda_\xi^{(n)}$, and $\varphi = \mu \circ \lambda_B^{(n)}$, then the sequences

$$Q_z^{B,n}\left(\frac{\varphi T}{n}, \frac{kT}{n}\right) - \check{Z}^B_{\theta^{(n)}_{\varphi \wedge k}}, \qquad k = 0, 1, \dots, n,$$

and

$$Q_z^{(n)}\left(\frac{\zeta T}{n}, \frac{kT}{n}\right) - \check{Z}^{\pi^\zeta, n}_{\zeta \wedge k}, \qquad k = 0, 1, \dots, n,$$

are obtained by means of the same functional from the sequences

$$B^*_{\theta_k^{(n)}} - B^*_{\theta_{k-1}^{(n)}}, \qquad k = 1, 2, \dots, \quad \text{and} \quad \left(\frac{T}{n}\right)^{1/2} \xi_k, \qquad k = 1, 2, \dots,$$

respectively. This together with (2.20) yields that

$$E^B \max_{0\le k\le n}\left(Q_z^{B,n}\left(\frac{\varphi T}{n},\frac{kT}{n}\right)-\check Z^B_{\theta_\varphi^{(n)}\wedge\theta_k^{(n)}}\right)^+=0. \tag{5.24}$$

In order to estimate the left-hand side of (2.22), we observe that

$$(R_z^B(\theta_\varphi^{(n)},t)-Z^B_{\theta_\varphi^{(n)}\wedge t})^+\le e^{rT}(Q_z^B(\theta_\varphi^{(n)},t)-\check Z^B_{\theta_\varphi^{(n)}\wedge t})^+ \tag{5.25}$$

and in view of (5.24),

$$\begin{aligned}(Q_z^B(\theta_\varphi^{(n)},t)-\check Z^B_{\theta_\varphi^{(n)}\wedge t})^+ &\le (Q_z^B(\theta_\varphi^{(n)},t)-Q_z^B(\theta_\varphi^{(n)},\theta_{\nu_t}^{(n)}))^+\\ &\quad+\left(Q_z^B(\theta_\varphi^{(n)},\theta_{\nu_t}^{(n)})-Q_z^{B,n}\left(\frac{\varphi T}{n},\frac{\nu_t T}{n}\right)\right)^+\\ &\quad+(\check Z^B_{\theta_\varphi^{(n)}\wedge\theta_{\nu_t}^{(n)}}-\check Z^B_{\theta_\varphi^{(n)}\wedge t})^+\mathbb{I}_{A_t},\end{aligned} \tag{5.26}$$

where we introduce the event $A_t=\{Q_z^B(\theta_\varphi^{(n)},t)\ge\check Z^B_{\theta_\varphi^{(n)}\wedge t}\}$ and, again, $\nu_t=\min\{k\in\mathbb{N}:\theta_k^{(n)}\ge t\}$. Note that since $\varphi\le n$,

$$\begin{aligned}Q_z^B(\theta_\varphi^{(n)},\theta_{\nu_t}^{(n)})&=Q_z^B(\theta_\varphi^{(n)},\theta_{\nu_t\wedge n}^{(n)}),\\ Q_z^{B,n}\left(\frac{\varphi T}{n},\frac{\nu_t T}{n}\right)&=Q_z^{B,n}\left(\frac{\varphi T}{n},(\nu_t\wedge n)\frac{T}{n}\right)\end{aligned} \tag{5.27}$$

and

$$Q_z^B(\theta_\varphi^{(n)}\wedge T,\theta_{\nu_t}^{(n)}\wedge T)=Q_z^B(\theta_\varphi^{(n)}\wedge T,\theta_{\nu_t\wedge n}^{(n)}\wedge T). \tag{5.28}$$

Taking into account (5.27) and (5.28), we obtain by the same estimates as in Lemma 3.4 that there exists $C^{(10)}>0$ such that, for all $z,n>0$,

$$\begin{aligned}&E^B\sup_{0\le t\le T}|Q_z^B(\theta_\varphi^{(n)},t)-Q_z^B(\theta_\varphi^{(n)}\wedge T,t)|\\ &\qquad\le C^{(10)}(F_0(z)+\Delta_0(z)+z+1)n^{-1/4}\end{aligned} \tag{5.29}$$

and

$$\begin{aligned}&E^B\sup_{0\le t\le T}|Q_z^B(\theta_\varphi^{(n)},\theta_{\nu_t\wedge n}^{(n)})-Q_z^B(\theta_\varphi^{(n)}\wedge T,\theta_{\nu_t\wedge n}^{(n)}\wedge T)|\\ &\qquad\le C^{(10)}(F_0(z)+\Delta_0(z)+z+1)n^{-1/4}.\end{aligned} \tag{5.30}$$

Similarly to the estimates in the proof of Lemma 3.6, we obtain that, for any $\varepsilon>0$, there exists $C_\varepsilon^{(11)}>0$ such that, for all $z,n>0$,

$$E^B\sup_{0\le t\le T}(Q_z^B(\theta_\varphi^{(n)}\wedge T,t)-Q_z^B(\theta_\varphi^{(n)}\wedge T,\theta_{\nu_t}^{(n)}\wedge T))^+$$

$$
\begin{aligned}
&\le E^B \sup_{0\le t\le T}(|e^{-\theta_\varphi^{(n)}\wedge t}-e^{-\theta^{(n)}_{\varphi\wedge\nu_t}\wedge T}|F_t(S^B(z)))\\
(5.31)\qquad &\quad + E^B\sup_{0\le t\le T}|F_t(S^B(z))-F_{\theta^{(n)}_{\nu_t}\wedge T}(S^B(z))|\\
&\le C_\varepsilon^{(11)} n^{\varepsilon-1/2}.
\end{aligned}
$$

Taking into account (5.27), we obtain by estimates of Lemmas 3.2 and 3.5 that there exists $C^{(12)}>0$ such that, for all $z,n>0$,

$$
\begin{aligned}
&E^B\sup_{0\le t\le T}\left(Q_z^B(\theta_\varphi^{(n)},\theta_{\nu_t}^{(n)})-Q_z^{B,n}\left(\frac{\varphi T}{n},\frac{\nu_t T}{n}\right)\right)^+\\
(5.32)\qquad &\quad\le C^{(12)}(F_0(z)+\Delta_0(z)+z+1)n^{-1/4}(\ln n)^{3/4}.
\end{aligned}
$$

By the definition, $\beta_s^\varphi=\beta^\varphi_{\theta^{(n)}_{\nu_s}}$ and $\gamma_s^\varphi=\gamma^\varphi_{\theta^{(n)}_{\nu_s}}$, and so

$$
\begin{aligned}
&(\check Z^B_{\theta_\varphi^{(n)}\wedge\theta_{\nu_t}^{(n)}}-\check Z^B_{\theta_\varphi^{(n)}\wedge t})^+\mathbb{I}_{A_t}\\
&\quad\le\gamma^\varphi_{\theta_\varphi^{(n)}\wedge t}(\check S^B_{\theta_\varphi^{(n)}\wedge\theta_{\nu_t}^{(n)}}(z)-\check S^B_{\theta_\varphi^{(n)}\wedge t}(z))^+\\
(5.33)\qquad &\quad\le\gamma^\varphi_{\theta_\varphi^{(n)}\wedge t}\check S^B_{\theta_\varphi^{(n)}\wedge t}(z)\mathbb{I}_{A_t}(\exp(\kappa(B^*_{\theta_\varphi^{(n)}\wedge\theta_{\nu_t}^{(n)}}-B^*_{\theta_\varphi^{(n)}\wedge t}))-1)^+\\
&\quad\le Q_z^B(\theta_\varphi^{(n)},t)\\
&\qquad\times\left(H_3^{(n)}+\mathbb{I}_{H_3^{(n)}>1}\left(1+e^{r\theta_n^{(n)}}\sup_{0\le s\le\theta_n^{(n)}}e^{\kappa B_s}\sup_{0\le s\le T}e^{-\kappa B_s}\right)\right),
\end{aligned}
$$

with $H_3^{(n)}$ defined before (4.43). Since, by (2.3) and (2.4),

$$
\sup_{0\le t\le T}Q_z^B(\theta_\varphi^{(n)},t)\le F_0(z)+\Delta_0(z)+L(T+2)\left(1+ze^{rT}\sup_{0\le t\le T}e^{\kappa B_s}\right),
$$

we derive from (4.8), (4.44), (4.45), (5.33) and the Hölder inequality that, for any $\varepsilon>0$, there exists $C_\varepsilon^{(13)}>0$ such that, for all $z,n>0$,

$$
\begin{aligned}
&E^B\sup_{0\le t\le T}(\check Z^B_{\theta_\varphi^{(n)}\wedge\theta_{\nu_t}^{(n)}}-\check Z^B_{\theta_\varphi^{(n)}\wedge t})^+\mathbb{I}_{A_t}\\
(5.34)\qquad &\quad\le C_\varepsilon^{(13)}(F_0(z)+\Delta_0(z)+z+1)n^{\varepsilon-1/2},
\end{aligned}
$$

which together with (5.25), (5.26) and (5.29)–(5.32) yields (2.22), completing the proof of Theorem 2.3. □

**6. Estimates á la Lamberton and Rogers.** In this section we derive a game option version of the approximation error estimates from [24] whose boundedness and smoothness assumptions do not permit to employ them even for standard options, but this explicit and simple method still has certain theoretical and pedagogical value. Let $\xi_1, \xi_2, \dots$ be i.i.d. random variables on a probability space $(\Omega^\xi, P^\xi)$ with $E^\xi \xi_1 = 0$, $E^\xi \xi_1^2 = 1$ and $E^\xi \xi_1^4 < \infty$. The latter ensures that if $\Theta$ is the Skorokhod embedding time of $\xi_1$ into the Brownian motion $B_t$ (i.e., a stopping time such that $\xi_1$ and $B_\Theta$ have the same distribution), then $\operatorname{Var}\Theta = E(\Theta - 1)^2 < \infty$ (see [4] or [24]). We will use the same notation as in Section 2 for $\Xi_k^{(n)}$ given by (2.17) and for the corresponding sets of stopping times $\mathcal{T}_{0T}^B$ and $\mathcal{T}_{0n}^\xi$ with respect to the Brownian filtration $\mathcal{F}_t^B$ with values in $[0,T]$ and with respect to the filtration $\mathcal{F}_k^\xi$ generated by $\xi_1, \xi_2, \dots$ with values in $\{1, \dots, n\}$, respectively. Let $g \ge f$ be continuous bounded functions on $[0,T] \times \mathbb{R}$ and

$$R(s,t,x) = g(s,x)\mathbb{I}_{s<t} + f(t,x)\mathbb{I}_{t\le s}.$$

Set

$$V = \inf_{\sigma \in \mathcal{T}_{0,T}^B} \sup_{\tau \in \mathcal{T}_{0,T}^B} E^B R(\sigma, \tau, B_{\sigma\wedge\tau}), \tag{6.1}$$

where, recall, $B_t, t \ge 0, B_0 = 0$ is the standard one-dimensional continuous in time Brownian motion, and

$$V^{(n)} = \inf_{\zeta \in \mathcal{T}_{0,n}^\xi} \sup_{\eta \in \mathcal{T}_{0,n}^\xi} E^\xi R\left(\frac{\zeta T}{n}, \frac{\eta T}{n}, \Xi_{\zeta\wedge\eta}^{(n)}\right). \tag{6.2}$$

THEOREM 6.1. *Let $f, g : \mathbb{R}^+ \times \mathbb{R} \to \mathbb{R}$ be bounded, continuous and having bounded and continuous derivatives $\frac{\partial f}{\partial t}$, $\frac{\partial g}{\partial t}$, $\frac{\partial^2 f}{\partial x^2}$ and $\frac{\partial^2 g}{\partial x^2}$ functions. Let*

$$Lh = \frac{\partial h}{\partial t} + \frac{1}{2}\frac{\partial^2 h}{\partial x^2}$$

*and $\rho = \sqrt{\operatorname{Var}\Theta}$. Then*

$$\begin{aligned} |V - V^{(n)}| \le{}& \frac{\rho T}{\sqrt{n}}\left(3\|Lf\|_\infty + 3\|Lg\|_\infty + \left\|\frac{\partial f}{\partial t}\right\|_\infty + \left\|\frac{\partial g}{\partial t}\right\|_\infty\right) \\ &+ \frac{T}{n}(\|Lf\|_\infty + \|Lg\|_\infty), \end{aligned} \tag{6.3}$$

*where $\|\cdot\|_\infty$ is the supremum norm on the whole $\mathbb{R}^+ \times \mathbb{R}$.*

PROOF. As in Section 2, we employ the Skorokhod embedding (see, e.g., [4]) which yields the existence of a nondecreasing sequence of stopping times $\theta_k^{(n)}$, $k = 1, 2, \dots$, $\theta_0^{(n)} = 0$ for the Brownian motion $B_t$ with its Brownian filtration $\mathcal{F}_t^B$, $t \ge 0$, such that $(\theta_{k+1}^{(n)} - \theta_k^{(n)}, B_{\theta_{k+1}^{(n)}} - B_{\theta_k^{(n)}})$ is independent of $\mathcal{F}^B_{\theta^{(n)}_k}$ and it has the same distribution as $(\frac{\Theta}{n}T, B_\Theta\sqrt{\frac{T}{n}})$. Let $\mathcal{T}^{B,n}$ be the set of integer valued stopping times with respect to the filtration $\{\mathcal{F}^B_{\theta^{(n)}_k}\}_{k\in\mathbb{N}}$ and the subset of these stopping times with values in $\{0,1,\dots,n\}$ we denote by $\mathcal{T}^{B,n}_{0,n}$. We claim that

$$V^{(n)} = V^{B,n} = \inf_{\zeta\in\mathcal{T}^{B,n}_{0,n}} \sup_{\eta\in\mathcal{T}^{B,n}_{0,n}} E^B R\left(\frac{\zeta T}{n}, \frac{\eta T}{n}, B_{\theta^{(n)}_\zeta\wedge\theta^{(n)}_\eta}\right). \tag{6.4}$$

Indeed, this result can be proved similarly to Lemma 3.1 employing the corresponding dynamical programming recursive formulas. Namely, we can write $V^{(n)} = V^{(n)}_{0,n}$, $V^{(n)}_{n,n} = f(T, \Xi^{(n)}_n)$ and for $k = 0, 1, \dots, n-1$,

$$V^{(n)}_{k,n} = \min\left(g\left(\frac{kT}{n}, \Xi^{(n)}_k\right), \max\left(f\left(\frac{kT}{n}, \Xi^{(n)}_k\right), E^\xi(V^{(n)}_{k+1,n}|\mathcal{F}^\xi_k)\right)\right)$$

and, on the other hand, $V^{B,n} = V^{B,n}_{0,n}$, $V^{B,n}_{n,n} = f(T, B_T)$ and for $k = 0, 1, \dots, n-1$,

$$V^{B,n}_{k,n} = \min\left(g\left(\frac{kT}{n}, B_{\theta^{(n)}_k}\right), \max\left(f\left(\frac{kT}{n}, B_{\theta^{(n)}_k}\right), E^B(V^{B,n}_{k+1,n}|\mathcal{F}^B_{\theta^{(n)}_k})\right)\right).$$

In the same way as in Lemma 3.1, we show by the backward induction that there exist a sequence $Q_k(x_1,\dots,x_k)$, $k = 1,\dots,n$, of measurable functions and a constant $Q_0$ such that $V^{(n)}_{0,n} = Q_0$, $V^{B,n}_{0,n} = Q_0$ and $V^{(n)}_{k,n} = Q_k((\frac{T}{n})^{1/2}\xi_1, \dots, (\frac{T}{n})^{1/2}\xi_k)$, $V^{B,n}_{kn} = Q_k(B_{\theta^{(n)}_1}, B_{\theta^{(n)}_2} - B_{\theta^{(n)}_1}, \dots, B_{\theta^{(n)}_k} - B_{\theta^{(n)}_{k-1}})$ for $k = 1,\dots,n$, and so (6.4) follows.

Next, set

$$V^{B,\theta,n} = \inf_{\zeta\in\mathcal{T}^{B,n}_{0,n}} \sup_{\eta\in\mathcal{T}^{B,n}_{0,n}} E^B R(\theta^{(n)}_\zeta, \theta^{(n)}_\eta, B_{\theta^{(n)}_\zeta\wedge\theta^{(n)}_\eta}). \tag{6.5}$$

Then

$$\begin{aligned} &|V^{B,n} - V^{B,\theta,n}| \\ &\quad\le \sup_{\zeta\in\mathcal{T}^{B,n}_{0,n}} \sup_{\eta\in\mathcal{T}^{B,n}_{0,n}} E^B\Big|R(\theta^{(n)}_\zeta, \theta^{(n)}_\eta, B_{\theta^{(n)}_\zeta\wedge\theta^{(n)}_\eta}) \\ &\qquad\qquad - R\left(\frac{\zeta T}{n}, \frac{\eta T}{n}, B_{\theta^{(n)}_\zeta\wedge\theta^{(n)}_\eta}\right)\Big| \\ &\quad\le \left(\left\|\frac{\partial f}{\partial t}\right\|_\infty + \left\|\frac{\partial g}{\partial t}\right\|_\infty\right) \sup_{\zeta\in\mathcal{T}^{B,n}_{0,n}} E^B\left|\theta^{(n)}_\zeta - \frac{\zeta T}{n}\right|. \end{aligned} \tag{6.6}$$

Since $\theta_k^{(n)} - \frac{kT}{n}$, $k = 0, 1, \dots, n$, is a martingale with respect to the filtration $\mathcal{F}^B_{\theta_k^{(n)}}$, $k \ge 0$, the sequence $|\theta_k^{(n)} - \frac{kT}{n}|$, $k = 0, 1, \dots, n$, is a submartingale which together with the Cauchy–Schwarz inequality yields

$$\begin{aligned} \sup_{\zeta \in \mathcal{T}^{B,n}_{0,n}} E^B \left|\theta_\zeta^{(n)} - \frac{\zeta T}{n}\right| &\le E^B |\theta_n^{(n)} - T| \\ &\le (E^B(\theta_n^{(n)} - T)^2)^{1/2} = \frac{\rho T}{\sqrt{n}}. \end{aligned} \tag{6.7}$$

Next, let $\mathcal{T}^{B,n}_T = \{\theta_\zeta^{(n)} \wedge T : \zeta \in \mathcal{T}^{B,n}_{0,n}\}$ and set

$$V^{B,n}_{0,T} = \inf_{\sigma \in \mathcal{T}^{B,n}_T} \sup_{\tau \in \mathcal{T}^{B,n}_T} E^B(e^{-r\sigma\wedge\tau} R(\sigma, \tau)). \tag{6.8}$$

By the Itô formula (see [11]), we arrive at the Dynkin formula, which gives

$$\begin{aligned} &|V^{B,\theta,n} - V^{B,n}_{0,T}| \\ &\quad \le \sup_{\zeta \in \mathcal{T}^{B,n}_{0,n}} \sup_{\eta \in \mathcal{T}^{B,n}_{0,n}} \Big| E^B \Big( R(\theta_\zeta^{(n)}, \theta_\eta^{(n)}, B_{\theta_\zeta^{(n)} \wedge \theta_\eta^{(n)}}) \\ &\qquad\qquad - R(\theta_\zeta^{(n)} \wedge T, \theta_\eta^{(n)} \wedge T, B_{\theta_\zeta^{(n)} \wedge \theta_\eta^{(n)} \wedge T}) \Big) \Big| \\ &\quad \le \sup_{\zeta \in \mathcal{T}^{B,n}_{0,n}} \left( \left| E^B \int_{\theta_\zeta^{(n)} \wedge T}^{\theta_\zeta^{(n)}} Lf(s, B_s)\, ds \right| + \left| E^B \int_{\theta_\zeta^{(n)} \wedge T}^{\theta_\zeta^{(n)}} Lg(s, B_s)\, ds \right| \right) \\ &\quad \le (\|Lf\|_\infty + \|Lg\|_\infty) E^B |\theta_n^{(n)} - T| \\ &\quad \le \rho \sqrt{\frac{T}{n}} (\|Lf\|_\infty + \|Lg\|_\infty). \end{aligned} \tag{6.9}$$

In order to obtain (6.3), it remains to show that

$$|V - V^{B,n}_{0,T}| \le \left( \frac{T}{n} + 2\frac{\rho T}{\sqrt{n}} \right) (\|Lf\|_\infty + \|Lg\|_\infty). \tag{6.10}$$

As in the proof of Lemma 3.6 in Section 4, for each $\sigma \in \mathcal{T}^B_{0,T}$, define $\nu_\sigma = \min\{k \in \mathbb{N} : \theta_k^{(n)} \ge \sigma\}$ and $\sigma^{(n)} = \theta_{\nu_\sigma}^{(n)} \wedge T$. Similarly to (4.58), we conclude that, for any $\delta > 0$, there exists $\sigma_\delta \in \mathcal{T}^B_{0,T}$ such that

$$\begin{aligned} V &\ge V^{B,n}_{0,T} - \delta - \sup_{\tau \in \mathcal{T}^{B,n}_T} E^B(R(\sigma_\delta^{(n)}, \tau, B_{\sigma_\delta^{(n)} \wedge \tau}) - R(\sigma_\delta, \tau, B_{\sigma_\delta \wedge \tau})) \\ &\ge V^{B,n}_{0,T} - \delta - \sup_{\sigma \in \mathcal{T}^B_{0,T}} |E^B(g(\sigma^{(n)}, B_{\sigma^{(n)}}) - g(\sigma, B_\sigma))| \end{aligned}$$

$$
\begin{aligned}
&- \sup_{\tau\in\mathcal{T}_{0,T}^B} |E^B(f(\tau^{(n)}, B_{\tau^{(n)}}) - g(\tau, B_\tau))| \\
&\geq V_{0,T}^{B,n} - \delta - (\|Lf\|_\infty + \|Lg\|_\infty) \sup_{\sigma\in\mathcal{T}_{0,T}^B} E^B|\sigma^{(n)} - \sigma|,
\end{aligned} \tag{6.11}
$$

where we used also the Itô formula and took into account in the same way as in (4.60) that $\sigma^{(n)} \geq \sigma$. By (4.62), for any $\sigma \in \mathcal{T}_{0,T}^B$,

$$
|\sigma^{(n)} - \sigma| \leq \max_{0\leq k\leq n-1} |\theta_{k+1}^{(n)} - \theta_k^{(n)}| + |T - \theta_n^{(n)}|,
$$

and so employing the estimate

$$
E \max_{0\leq k\leq n-1} Z_k \leq EZ_0 + \sqrt{n\,\mathrm{Var}\,Z_0}
$$

from Lemma 3.5 of [24] to the i.i.d. random variables $Z_k = \theta_{k+1}^{(n)} - \theta_k^{(n)}$, we obtain from (4.1) and (6.7) that

$$
E^B|\sigma^{(n)} - \sigma| \leq \frac{T}{n} + 2\frac{T\rho}{\sqrt{n}}. \tag{6.12}
$$

Taking into account that $\delta$ in (6.11) can be taken arbitrarily small, we obtain from (6.11) and (6.12) that

$$
V - V_{0,T}^{B,n} \geq -\left(\frac{T}{n} + 2\frac{T\rho}{\sqrt{n}}\right)(\|Lf\|_\infty + \|Lg\|_\infty). \tag{6.13}
$$

By [23], we can write also that

$$
\begin{aligned}
V &= \sup_{\tau\in\mathcal{T}_{0,T}^B} \inf_{\sigma\in\mathcal{T}_{0,T}^B} E^B R(\sigma, \tau, B_{\sigma\wedge\tau}) \\
&\leq V_{0,T}^{B,n} + \delta + \sup_{\sigma\in\mathcal{T}_T^{B,n}} (R(\sigma, \tau_\delta, B_{\sigma\wedge\tau_\delta}) - R(\sigma, \tau_\delta^{(n)}, B_{\sigma\wedge\tau_\delta^{(n)}}))
\end{aligned}
$$

and similarly to the above, we obtain that

$$
V - V_{0,T}^{B,n} \leq \left(\frac{T}{n} + 2\frac{T\rho}{\sqrt{n}}\right)(\|Lf\|_\infty + \|Lg\|_\infty),
$$

which together with (6.13) gives (6.10) and completes the proof of Theorem 6.1. □

**Acknowledgments.** I would like to thank the anonymous referees for carefully reading this paper and indicating to me a number of errors and misprints in its first version.

Institute of Mathematics
Hebrew University
Jerusalem 91904
Israel
E-mail: kifer@math.huji.ac.il